\documentclass[12pt]{article}

\usepackage{amssymb}
\usepackage{amsmath}
\usepackage{cite}

\begin{document}

\renewcommand{\citeleft}{{\rm [}}
\renewcommand{\citeright}{{\rm ]}}
\renewcommand{\citepunct}{{\rm,\ }}
\renewcommand{\citemid}{{\rm,\ }}

\newcounter{abschnitt}
\newtheorem{satz}{Theorem}
\newtheorem{theorem}{Theorem}[abschnitt]
\newtheorem{koro}[theorem]{Corollary}
\newtheorem{prop}[theorem]{Proposition}
\newtheorem{lem}[theorem]{Lemma}
\newtheorem{conj}[theorem]{Conjecture}

\newcounter{saveeqn}
\newcommand{\alpheqn}{\setcounter{saveeqn}{\value{abschnitt}}
\renewcommand{\theequation}{\mbox{\arabic{saveeqn}.\arabic{equation}}}}
\newcommand{\reseteqn}{\setcounter{equation}{0}
\renewcommand{\theequation}{\arabic{equation}}}

\hyphenation{convex} \hyphenation{bodies}

\hyphenpenalty=9000

\sloppy

\phantom{a}

\vspace{-1.7cm}

\begin{center}
 \begin{Large} {\bf Even Minkowski Valuations} \\[0.6cm] \end{Large}

\begin{large} Franz E. Schuster and Thomas Wannerer \end{large}
\end{center}

\vspace{-1cm}

\begin{quote}
\footnotesize{ \vskip 1truecm\noindent {\bf Abstract.} A new
integral representation of smooth translation invariant and rotation equivariant even Minkowski valuations
is established. Explicit formulas
relating previously obtained descriptions of such valuations with
the new more accessible one are also derived. Moreover, the
action of Alesker's Hard Lefschetz operators on these Minkowski valuations is explored in detail.}
\end{quote}

\vspace{0.6cm}

\centerline{\large{\bf{ \setcounter{abschnitt}{1}
\arabic{abschnitt}. Introduction}}} \alpheqn

\vspace{0.6cm}

A \emph{Minkowski valuation} is a map $\Phi: \mathcal{K}^n \rightarrow \mathcal{K}^n$ defined on the set $\mathcal{K}^n$
of convex bodies (compact convex sets) in $\mathbb{R}^n$ such that
\[\Phi K + \Phi L = \Phi(K \cup L) + \Phi(K \cap L), \]
whenever $K \cup L \in \mathcal{K}^n$ and addition on $\mathcal{K}^n$ is the usual Minkowski addition.
While a number of classical (reverse) affine isoperimetric inequalities
involve \linebreak well known Minkowski valuations, like the projection and
difference body maps, the underlying reason for the special role of these maps
has only been revealed recently, when they
were characterized as the unique Minkowski valuations which
intertwine linear transformations. This line of research can be
traced back to two seminal articles by Ludwig
\textbf{\cite{ludwig02, Ludwig:Minkowski}} and has become the
focus of increased interest in recent years, see
\textbf{\cite{abardber11, haberl10, haberl09, hablud06, haberl11,
Ludwig06, Ludwig10a, Ludwig11, SchuWan11, wannerer10}}.

Due to these recent characterization results a more and more complete picture on linearly intertwining
Minkowski valuations could be developed. However, the theory of Minkowski
valuations which are merely translation invariant and $\mathrm{SO}(n)$ equivariant is still in its
infancy. Here, Schneider \textbf{\cite{schneider74}} and Kiderlen
\textbf{\cite{kiderlen05}} have obtained first
classification theorems of Minkowski valuations homogeneous of degree one. For
{\it even} valuations their results were subsequently
generalized by the first author \textbf{\cite{Schu09}}:
It was shown that smooth translation invariant and $\mathrm{SO}(n)$ equivariant even Minkowski
valuations of an {\it arbitrary} degree of homogeneity are
generated by convolutions between projection functions
and {\it spherical Crofton measures}. Unfortunately, the rather
complicated technical nature of this general result also made it
difficult to work with. In particular, uniqueness as well as the
problem of finding necessary and sufficient conditions for a
measure to be a spherical Crofton measure of a Minkowski valuation remained
(essentially) open.

In this article we first show that the spherical Crofton measure of an even
Minkowski valuation is uniquely
determined. This is done by exploiting tools from harmonic
analysis and a different description of such valuations based on
an embedding of even translation invariant real valued valuations
in continuous functions on the Grassmannian by Klain
\textbf{\cite{klain00}}. Our approach has the added advantage
that it helps to simplify a previously known necessary condition
for a measure to be the spherical Crofton
 measure of a Minkowski valuation: We show that its
spherical cosine transform has to be the support function of a
convex body. As our main result we then derive a new more
elementary integral representation of smooth translation invariant and $\mathrm{SO}(n)$ equivariant even Minkowski
valuations. We also give explicit
formulas, involving Radon transforms on Grassmannians, which
relate the previously known descriptions of such valuations with
our new more accessible one.

By recent results of Bernig \textbf{\cite{Bernig07}} and
Parapatits with the first author \textbf{\cite{parapschu}}
Alesker's Hard Lefschetz operators (originally defined only for
translation invariant real valued valuations; see
\textbf{\cite{Alesker03, Alesker04, alesker10, bernigfu06}}) can
be extended to translation invariant and $\mathrm{SO}(n)$ equivariant Minkowski valuations. In the final section of this article we study these
operators in terms of their action on the different possible
representations of such Minkowski valuations.

\vspace{0.7cm}

\centerline{\large{\bf{ \setcounter{abschnitt}{2}
\arabic{abschnitt}. Statement of principal results}}}

\reseteqn \alpheqn

\vspace{0.4cm}

We endow the set $\mathcal{K}^n$ of convex bodies in $\mathbb{R}^n$ with the
Hausdorff metric and assume throughout that $n \geq 3$. All measures in this article are
\emph{signed} finite Borel measures. A convex body $K$ is uniquely determined by its
support function $h(K,u)=\max \{u \cdot x: x \in K\}$ for $u \in
S^{n-1}$. For $i \in \{1, \ldots, n - 1\}$, the $i$th projection
function $\mathrm{vol}_i(K|\,\cdot\,)$ of $K \in \mathcal{K}^n$
is the continuous function on the Grassmannian
$\mathrm{Gr}_{i,n}$ of $i$-dimensional subspaces of
$\mathbb{R}^n$ defined such that $\mathrm{vol}_i(K|E)$, for $E
\in \mathrm{Gr}_{i,n}$, is the $i$-dimensional volume of the
orthogonal projection of $K$ onto $E$.

A function $\phi: \mathcal{K}^n \rightarrow \mathcal{A}$ with values in an Abelian semigroup $\mathcal{A}$ is called a {\it
valuation}, or {\it additive}, if
\[\phi(K) + \phi(L) = \phi(K \cup L) + \phi(K \cap L)   \]
whenever $K \cup L$ is convex. With its origins in Dehn's
solution of Hilbert's Third Problem, the notion of {\it scalar}
valued valuations has long played a central role
in convex, discrete, and integral geometry (see, e.g.,
\textbf{\cite{Klain:Rota}} for more information). Minkowski valuations are of newer vintage.

Recent classifications of Minkowski valuations (see
\textbf{\cite{haberl11, ludwig02, Ludwig:Minkowski, SchuWan11,
wannerer10}}) showed that only a small number of such operators,
like the projection and the difference body map (see Section 4 for
their definitions), intertwine \emph{affine} transformations. In
this article, we study the much larger class of continuous
Minkowski valuations which are \emph{translation invariant} and
\emph{$\mathrm{SO}(n)$ equivariant}.

A map $\Phi$ from $\mathcal{K}^n$ to $\mathcal{K}^n$ (or
$\mathbb{R}$) is said to have {\it degree} $i$ if $\Phi(\lambda
K) = \lambda^i\Phi K$ \linebreak  for $K \in \mathcal{K}^n$ and $\lambda >
0$. In the case of Minkowski valuations which are translation invariant and
$\mathrm{SO}(n)$ equivariant of degree $i$, Kiderlen \textbf{\cite{kiderlen05}} \linebreak (for $i = 1$;
building on previous results by Schneider
\textbf{\cite{schneider74}}) and the first author
\textbf{\cite{Schu06a}} (for $i = n - 1$) were the first to
obtain representations of these maps by spherical convolution
operators. The following extension of their results to all
remaining (non-trivial) degree cases (by a result of McMullen
\textbf{\cite{McMullen77}}, only integer degrees $0 \leq i \leq
n$ can occur) was established in \textbf{\cite{Schu09}} for
\emph{even} Minkowski valuations, that is, \ $\Phi(-K)=\Phi K$ for
$K \in \mathcal{K}^n$:

\begin{satz}[\!\!\cite{Schu09}] \label{duke} Let $\Phi_i: \mathcal{K}^n \rightarrow
\mathcal{K}^n$ be a smooth translation invariant and \linebreak
$\mathrm{SO}(n)$ equivariant Minkowski valuation of degree $i \in
\{1, \ldots, n - 1\}$. If $\Phi_i$ is even, then there exists a
smooth $\mathrm{O}(i) \times \mathrm{O}(n - i)$ invariant measure
$\mu$ on $S^{n-1}$ such that for every $K \in \mathcal{K}^n$,
\begin{equation} \label{crofton}
h(\Phi_iK,\cdot) = \mathrm{vol}_i(K|\,\cdot\,) \ast \mu.
\end{equation}
\end{satz}

We note that Theorem \ref{duke} was stated in a different, but
equivalent, form in \textbf{\cite{Schu09}} (cf.\ the discussion
in Section 4 and the Appendix). The convolution in (\ref{crofton})
is induced from the group $\mathrm{SO}(n)$ by identifying
$S^{n-1}$ and $\mathrm{Gr}_{i,n}$ with the homogeneous spaces
$\mathrm{SO}(n)/\mathrm{SO}(n - 1)$ and
$\mathrm{SO}(n)/\mathrm{S}(\mathrm{O}(i) \times \mathrm{O}(n-i))$,
respectively (see Section 3 for details).

The notion of {\it smooth} translation invariant Minkowski
valuations which are $\mathrm{SO}(n)$ equivariant was introduced
in \textbf{\cite{Schu09}} (extending the definition of smooth
scalar valued valuations by Alesker \textbf{\cite{Alesker03}};
see Section 4). Moreover, it was shown in \textbf{\cite{Schu09}}
that every translation invariant and
$\mathrm{SO}(n)$ equivariant Minkowski valuation which is
continuous and even can be approximated uniformly on compact
subsets of $\mathcal{K}^n$ by smooth ones.

The invariant signed measure appearing in Theorem \ref{duke} is
essentially a Crofton measure of a real valued valuation
associated with the given even Minkowski valuation (see Section
4). This motivates the following definition.

\vspace{0.3cm}

\noindent {\bf Definition} \emph{Let $\Phi_i: \mathcal{K}^n
\rightarrow \mathcal{K}^n$ be an even Minkowski valuation of
degree $i \in \{1, \ldots, n - 1\}$. We call an $\mathrm{O}(i)
\times \mathrm{O}(n - i)$ invariant measure $\mu$ on $S^{n-1}$ a
\emph{spherical Crofton measure} for $\Phi_i$ if (\ref{crofton}) holds for
every $K \in \mathcal{K}^n$.}

\vspace{0.4cm}

Note that an even Minkowski valuation which admits a spherical Crofton
measure is continuous, translation invariant, and
$\mathrm{SO}(n)$ equivariant.

While Theorem \ref{duke} generalized the
previously known representation results \linebreak of translation invariant and
$\mathrm{SO}(n)$ equivariant even Minkowski
valuations, basic questions concerning
spherical Crofton measures, such as their uniqueness, remained open. With our first
result we answer some of them:

\begin{satz} \label{main1} If $\Phi_i: \mathcal{K}^n
\rightarrow \mathcal{K}^n$ is an even Minkowski valuation of
degree $i \in \{1, \ldots, n - 1\}$ which admits a spherical Crofton
measure $\mu$, then $\Phi_i$ determines $\mu$ uniquely. Moreover,
there exists an $\mathrm{O}(i) \times \mathrm{O}(n-i)$ invariant
convex body $L \in \mathcal{K}^n$ such that
\begin{equation} \label{kb}
h(L,u) = \int_{S^{n-1}} |u \cdot v|\, d\mu(v), \qquad u \in S^{n-1}.
\end{equation}
\end{satz}

The proof of Theorem~\ref{main1} relies on tools from harmonic
analysis and Klain's \textbf{\cite{klain00}} embedding of even
translation invariant real valued valuations in continuous
functions on the Grassmannian (see Section 4). In fact, the
support function in (\ref{kb}) is closely related to the Klain
function of the real valued valuation associated with $\Phi_i$.

Note that if a spherical Crofton measure $\mu$ exists for an even Minkowski
valuation $\Phi_i$, then, by its uniqueness and the injectivity of
the spherical cosine transform, also the body $L$ from (\ref{kb})
is uniquely determined by $\Phi_i$. In \textbf{\cite{Schu09}} it
was shown that, while not every \emph{continuous} translation invariant and $\mathrm{SO}(n)$ equivariant even Minkowski
valuation $\Phi_i$ admits a spherical
Crofton measure, it is always possible to associate a convex body
with $\Phi_i$ that uniquely determines the Minkowski valuation. We
call this body the {\it Klain body} of $\Phi_i$. \linebreak If $\Phi_i$
admits a spherical Crofton measure $\mu$, then its Klain body is (up to a
factor) the body $L$ given by (\ref{kb}). However, we recall in
Section 4 a much simpler way to determine this body from the
given Minkowski valuation.

\vspace{0.2cm}

Using Theorems \ref{duke} and \ref{main1}, we establish in
Section 5 the following new and much more elementary integral
representation for smooth translation invariant and
$\mathrm{SO}(n)$ equivariant even Minkowski valuations:

\begin{satz} \label{main2} Let $\Phi_i: \mathcal{K}^n \rightarrow
\mathcal{K}^n$ be a smooth translation invariant and
$\mathrm{SO}(n)$ equivariant Minkowski valuation of degree $i \in
\{1, \ldots, n - 1\}$. If $\Phi_i$ is even, then there exists a
unique even $g \in C^{\infty}((-1,1)) \cap C([-1,1])$ such that
for every $K \in \mathcal{K}^n$,
\begin{equation} \label{area}
h(\Phi_iK,u) = \int_{S^{n-1}} g(u \cdot v)\,dS_i(K,v), \qquad u
\in S^{n-1}.
\end{equation}
\end{satz}

Here, the measures $S_i(K,\cdot), 1 \leq i \leq n - 1$, are the
area measures of order $i$ of $K \in \mathcal{K}^n$ (see Section
4 for their definition).

Now, given three different ways to represent a smooth translation invariant and
$\mathrm{SO}(n)$ equivariant even Minkowski valuation (namely,
spherical Crofton measures, Klain bodies, and the {\it generating
functions} on $[-1,1]$ from Theorem \ref{main2}), a natural question to ask is how to convert one into another. A
partial answer is already provided by Theorem \ref{main1}: The
Klain body can be obtained from the spherical Crofton measure by the
spherical cosine transform. In Section 5, we derive explicit
formulas involving Radon transforms on Grassmannians for all the
remaining cases.

\vspace{0.2cm}

In the final section we explore the action of Alesker's Hard
Lefschetz operators on translation invariant and
$\mathrm{SO}(n)$ equivariant even Minkowski valuations. These operators, one being a derivation and
the other an integration operator, have been introduced
only recently and have since played an important role in what is
now called algebraic integral geometry (see, e.g.,
\textbf{\cite{Alesker03, bernigfu10}}). The Hard Lefschetz
Theorem for them (see \textbf{\cite{Alesker03, Alesker04,
alesker10, Bernig07b}}) is a fundamental theorem in the theory of
translation invariant {\it scalar valued} valuations. The authors
feel that the Hard Lefschetz operators will play a similarly
important role in the theory  of convex body valued valuations,
in particular, in connection with geometric inequalities (see
\textbf{\cite{BPSW2014, parapschu}} for first results in this direction). A
critical tool in our investigations is a Fourier type transform,
introduced for translation invariant scalar valued valuations by
Alesker \textbf{\cite{alesker10}}, that connects the two Hard
Lefschetz operators.

\vspace{1cm}

\centerline{\large{\bf{ \setcounter{abschnitt}{3}
\arabic{abschnitt}. Intertwining transforms on Grassmannians}}}

\reseteqn \alpheqn

\vspace{0.6cm}

In the following we recall the notion of convolution of functions
on the homogeneous spaces $\mathrm{SO}(n)/\mathrm{SO}(n - 1)$ and
$\mathrm{SO}(n)/\mathrm{S}(\mathrm{O}(i) \times \mathrm{O}(n-i))$
as well as basic facts about spherical functions on
Grass\-mannians. At the end of this section we establish two
critical auxiliary results which are needed in the proofs of our
main results. As a reference for this section we recommend the
book \textbf{\cite{takeuchi}} by Takeuchi and the article
\textbf{\cite{grinbergzhang99}} by Grinberg and Zhang.

In order to simplify the exposition we first consider a general
compact Lie group $G$ and a closed subgroup $H$ of $G$ (although
we only need the cases where $G = \mathrm{SO}(n)$ and $H$ is
either $\mathrm{SO}(n - 1)$ or $\mathrm{S}(\mathrm{O}(i) \times
\mathrm{O}(n-i))$).

\pagebreak

Let $C(G)$ denote the space of continuous functions on $G$. For $f, g \in C(G)$, the convolution $f
\ast g \in C(G)$ is defined by
\[(f \ast g)(\eta)=\int_G f(\eta \vartheta^{-1})g(\vartheta)\,d\vartheta
=\int_Gf(\vartheta)g(\vartheta^{-1}\eta)\,d\vartheta,\] where
integration is with respect to the Haar probability measure on
$G$.

For $f \in C(G)$ and a measure $\mu$ on $G$, the convolutions $f \ast \mu \in C(G)$ and $\mu \ast
f \in C(G)$ are defined by
\begin{equation} \label{convmeas}
(f \ast \mu)(\eta)=\int_G f(\eta
\vartheta^{-1})\,d\mu(\vartheta), \qquad (\mu \ast f)(\eta)=\int_G
f(\vartheta^{-1}\eta)\,d\mu(\vartheta).
\end{equation}
This definition continuously extends the convolution of
functions. Moreover, it follows from (\ref{convmeas}) that $f \ast
\mu$ and $\mu \ast f$ are smooth if $f \in C^{\infty}(G)$.

For $\vartheta \in G$, we denote the left and right translations
of $f \in C(G)$ by
\begin{equation*}
(l_{\vartheta}f)(\eta)=f(\vartheta^{-1}\eta), \qquad
(r_{\vartheta}f)(\eta)=f(\eta\vartheta).
\end{equation*}
For a measure $\mu$ on $G$, we define $l_{\vartheta}\mu$ and
$r_{\vartheta}\mu$ as the image measures of $\mu$ under
the left and right multiplication by $\vartheta$ in $G$,
respectively. We will also often write $\vartheta f := l_{\vartheta}f$ and
$\vartheta\mu := l_{\vartheta}\mu$ for the left translation of $f
\in C(G)$ or a measure $\mu$ on $G$, respectively.

We emphasize that, if $f \in C(G)$ and $\mu$ is a measure on $G$, then, by
(\ref{convmeas}),
\begin{equation} \label{roteq}
(l_{\vartheta} f) \ast \mu = l_{\vartheta}(f \ast \mu) \qquad
\mbox{and} \qquad \mu \ast (r_{\vartheta} f) = r_{\vartheta}(\mu
\ast f).
\end{equation}
for every $\vartheta \in
G$. Thus, the convolution from the right gives rise to operators on
$C(G)$ which intertwine left translations and the convolution
from the left gives rise to operators which intertwine right
translations.

For $f \in C(G)$, the function $\widehat{f} \in C(G)$ is defined
by
\begin{equation*}
\widehat{f}(\vartheta)=f(\vartheta^{-1}).
\end{equation*}
For a measure $\mu$ on $G$, we define the measure $\widehat{\mu}$
by
\begin{equation*}
\int_G f(\vartheta)\,d\widehat{\mu}(\vartheta) = \int_G
f(\vartheta^{-1})\,d\mu(\vartheta), \qquad f \in C(G).
\end{equation*}

From (\ref{convmeas}) it follows that for $f, g \in C(G)$ and a
measure $\sigma$ on $G$,
\begin{equation} \label{adjoint}
\int_G f(\vartheta)(g \ast \sigma)(\vartheta)\,d\vartheta =
\int_G (f \ast
\widehat{\sigma})(\vartheta)g(\vartheta)\,d\vartheta.
\end{equation}
Therefore, it is consistent to define the convolution $\mu \ast
\sigma$ of two measures $\mu, \sigma$ on $G$ by
\begin{equation*}
\int_G f(\vartheta)\,d(\mu \ast \sigma)(\vartheta)=\int_G (f\ast
\widehat{\sigma})(\vartheta)\,d\mu(\vartheta)=\int_G
(\widehat{\mu} \ast f)(\vartheta)\,d\sigma(\vartheta)
\end{equation*}
for every $f \in C(G)$. The convolution of functions and measures
on $G$, thus defined, is easily seen to be associative but in
general not commutative. If $\mu, \sigma$ are measures on $G$,
then
\begin{equation} \label{commute}
\widehat{\mu \ast \sigma}=\widehat{\sigma} \ast \widehat{\mu}.
\end{equation}

We now define convolutions between
functions and measures on the \linebreak homogeneous spaces
\begin{equation} \label{ident17}
S^{n-1}=\mathrm{SO}(n)/\mathrm{SO}(n-1) \quad \mbox{and} \quad \mathrm{Gr}_{i,n}=\mathrm{SO}(n)/\mathrm{S}(\mathrm{O}(i) \times \mathrm{O}(n-i))
\end{equation}
by identifying the space $C(G/H)$ of continuous functions on the homogeneous space $G/H$ with the closed subspace of $C(G)$ of all functions which are
right $H$-invariant, that is, $r_{\vartheta}f = f$ for every $\vartheta \in H$. Similarly, we identify measures on $G/H$ with right $H$-invariant measures on
$G$. For a detailed description of the one-to-one correspondence between functions and measures on $G/H$ and right $H$-invariant
functions and measures on $G$ we refer to \textbf{\cite{grinbergzhang99}} or \textbf{\cite{Schu09}}.

If, for example, $f \in
C(G/H)$ and $\mu$ is a measure on $G/H$, then, by (\ref{roteq}), the
convolution $f \ast \mu$ satisfies
\begin{equation} \label{convhom}
r_{\vartheta}(f \ast \mu)=f \ast(r_{\vartheta} \mu) = f \ast \mu
\end{equation}
for every $\vartheta \in G$. Thus $f \ast \mu$ can again be identified with
a continuous function on $G/H$. In the same way we can define convolutions between functions and measures on different
homogeneous spaces: Let $H_1, H_2$ be two closed subgroups of
$G$. If, say, $f \in C(G/H_1)$ and $g \in
C(G/H_2)$, then, by (\ref{roteq}), $f \ast g$ defines a
continuous right $H_2$-invariant function on $G$ and,
thus, can be identified with a continuous function on
$G/H_2$.

Let $\pi: G \rightarrow G/H$ be the
canonical projection and write $\pi(\vartheta)= \bar{\vartheta}$.
If $e \in G$ denotes the identity element, then $H$ is the
stabilizer in $G$ of $\bar{e} \in G/H$
and we have $\bar{\vartheta}=\vartheta \bar{e}$ for every
$\vartheta \in G$.

If $\delta_{\bar{e}}$ denotes the Dirac measure on $G/H$, then it is not difficult to show that for
$f \in C(G)$,
\begin{equation} \label{dirac1}
f \ast \delta_{\bar{e}} = \int_H r_{\vartheta}f\,d\vartheta.
\end{equation}
Thus, $f \ast \delta_{\bar{e}}$ is right $H$-invariant for every $f \in C(G)$ and $\delta_{\bar{e}}$ is the unique rightneutral
element for the convolution of functions and measures on $G/H$. We also note that
\begin{equation} \label{dirac2}
\delta_{\bar{e}} \ast f  =\int_{H}\vartheta f\,d\vartheta
\end{equation}
is left $H$-invariant for every $f \in C(G)$.

We call a left $H$-invariant function or measure on $G/H$ (or, equivalently, an $H$-biinvariant function or measure on $G$) {\it zonal}.
Zonal functions (and measures) on
$G/H$ play an essential role with respect to convolutions:
If $f, g \in C(G/H)$, then, by (\ref{dirac1})
and (\ref{dirac2}),
\begin{equation} \label{zonsuff}
f \ast g = (f \ast \delta_{\bar{e}}) \ast g = f \ast
(\delta_{\bar{e}} \ast g).
\end{equation}
Hence, for convolutions from the right on
$G/H$, it is sufficient to consider zonal functions
and measures. Note that if $f \in C(G)$ is $H$-biinvariant (or, equivalently, $f \in C(G/H)$ is zonal), then the function
$\widehat{f} \in C(G)$ is also $H$-biinvariant and, thus, zonal.

A different consequence of the identifications (\ref{ident17}) which we will use frequently is the following: If $h \in C(S^{n-1})$ is $\mathrm{S}(\mathrm{O}(i) \times \mathrm{O}(n-i))$ invariant,
then $\widehat{h} \in C(\mathrm{Gr}_{i,n})$ is $\mathrm{SO}(n - 1)$ invariant and, vice versa,
if $f \in C(\mathrm{Gr}_{i,n})$ is $\mathrm{SO}(n - 1)$ invariant, then $\widehat{f} \in C(S^{n-1})$ is $\mathrm{S}(\mathrm{O}(i) \times \mathrm{O}(n-i))$ invariant.

The following examples of convolution transforms will play a critical role in the proof of our main theorems.

\vspace{0.3cm}

\noindent {\bf Examples:}
\begin{enumerate}
\item[(a)] Suppose that $1 \leq i \neq j \leq n - 1$ and let $F \in \mathrm{Gr}_{j,n}$. We denote by $\mathrm{Gr}_{i,n}^F$ the submani\-fold of
$\mathrm{Gr}_{i,n}$ which comprises of all $E \in \mathrm{Gr}_{i,n}$ that contain (respectively, are contained in) $F$. The
{\it Radon transform} \linebreak $R_{i,j}: L^2(\mathrm{Gr}_{i,n}) \rightarrow L^2(\mathrm{Gr}_{j,n})$ is defined by
\[(R_{i,j}f)(F) = \int_{\mathrm{Gr}_{i,n}^F}f(E)\,d\nu_i^F(E),   \]
where $\nu_i^F$ is the unique invariant probability measure on $\mathrm{Gr}_{i,n}^F$.
It is well known that the Radon transform is a continuous linear operator and that $R_{j,i}$ is the
adjoint of $R_{i,j}$, that is,
\[\int_{\mathrm{Gr}_{j,n}} (R_{i,j}f)(F)g(F)\,dF = \int_{\mathrm{Gr}_{i,n}}f(E)(R_{j,i}g)(E)\,dE   \]
for every $f \in L^2(\mathrm{Gr}_{i,n})$ and $g \in L^2(\mathrm{Gr}_{j,n})$.

\pagebreak

Using the last observation, one can extend the Radon transform to the space of measures on $\mathrm{Gr}_{i,n}$ by
\[\int_{\mathrm{Gr}_{j,n}}f(F)\,d(R_{i,j}\mu)(F) = \int_{\mathrm{Gr}_{i,n}} (R_{j,i}f)(E)\,d\mu(E)\]
for $f \in C(\mathrm{Gr}_{j,n})$. For $g \in L^2(\mathrm{Gr}_{i,n})$, we write $g^{\bot}$ for the function in
$L^2(\mathrm{Gr}_{n-i,n})$ defined by $g^{\bot}(E) = g(E^{\bot})$. Then
\begin{equation} \label{rijbot}
(R_{i,j}f)^{\bot} = R_{n - i,n - j}f^{\bot}
\end{equation}
and, for $1 \leq i < j < k \leq n - 1$,  we have
\[R_{i,k} = R_{j,k} \circ R_{i,j} \qquad \mbox{and} \qquad R_{k,i} = R_{j,i} \circ R_{k,j}.   \]
If $1 \leq i < j \leq n - 1$ and $\lambda_{i,j}$ denotes the probability measure on $\mathrm{Gr}_{j,n}$ which is uniformly concentrated on the submanifold
\[\{\vartheta \bar{e} \in \mathrm{Gr}_{j,n}: \vartheta \in \mathrm{S}(\mathrm{O}(i) \times \mathrm{O}(n - i))\},   \]
then (see, e.g., \textbf{\cite{grinbergzhang99}})
\begin{equation} \label{radconv}
R_{i,j}f = f \ast \lambda_{i,j} \qquad \mbox{and} \qquad R_{j,i}g = g \ast \widehat{\lambda}_{i,j}
\end{equation}
for every $f \in L^2(\mathrm{Gr}_{i,n})$ and $g \in
L^2(\mathrm{Gr}_{j,n})$. In particular, the Radon transform
intertwines the natural group action (by left translation) of
$\mathrm{SO}(n)$ and maps smooth functions to smooth ones, that
is,
\[R_{i,j}: C^{\infty}(\mathrm{Gr}_{i,n}) \rightarrow C^{\infty}(\mathrm{Gr}_{j,n}).   \]
\item[(b)] For two subspaces $E, F \in \mathrm{Gr}_{i,n}$, where $1 \leq i \leq n - 1$, the cosine of the angle between $E$ and $F$ is
defined by
\[|\cos(E,F)| = \mathrm{vol}_i(Q|E),   \]
where $Q$ is an arbitrary subset of $F$ with $\mathrm{vol}_i(Q)=1$. (This definition is independent of the choice of $Q \subseteq F$.)
The continuous linear operator $C_i: L^2(\mathrm{Gr}_{i,n}) \rightarrow L^2(\mathrm{Gr}_{i,n})$ defined by
\[(C_if)(F) = \int_{\mathrm{Gr}_{i,n}} |\cos(E,F)|f(E)\,dE  \]
is called the {\it cosine transform}. Here integration is with respect to the Haar probability measure on $\mathrm{Gr}_{i,n}$.

\pagebreak

It is easy to see that also the cosine transform is a continuous
linear operator and that it is self-adjoint, that is,
\begin{equation} \label{ciselfad}
\int_{\mathrm{Gr}_{i,n}} (C_if)(E)g(E)\,dE =
\int_{\mathrm{Gr}_{i,n}}f(E)(C_ig)(E)\,dE
\end{equation}
for all $g, f \in L^2(\mathrm{Gr}_{i,n})$. Based on
(\ref{ciselfad}), one defines the cosine transform of a measure
$\mu$ on $\mathrm{Gr}_{i,n}$ by
\[\int_{\mathrm{Gr}_{i,n}}f(E)\,d(C_i\mu)(E) = \int_{\mathrm{Gr}_{i,n}} (C_if)(E)\,d\mu(E)\]
for $f \in C(\mathrm{Gr}_{i,n})$. For all $f \in
L^2(\mathrm{Gr}_{i,n})$, we have
\begin{equation} \label{cifbot}
(C_if)^{\bot} = C_{n-i}f^{\bot}
\end{equation}
and in \textbf{\cite{goodeyzhang}} it was shown that for all $1
\leq i \neq j \leq n - 1$,
\begin{equation} \label{rijci}
R_{i,j} \circ C_i = \frac{i!(n-i)!\kappa_i \kappa_{n-i}}{j!(n-j)!\kappa_j \kappa_{n-j}}\,C_j \circ R_{i,j},
\end{equation}
where $\kappa_i$ is the $i$-dimensional volume of the
$i$-dimensional Euclidean unit ball. It is not difficult to show that
\begin{equation} \label{cosconv}
C_i f = f \ast |\cos(\bar{e},\cdot)|
\end{equation}
for every $f \in L^2(\mathrm{Gr}_{i,n})$. Thus, also the cosine
transform intertwines the action of $\mathrm{SO}(n)$ and it
follows that
\[C_i: C^{\infty}(\mathrm{Gr}_{i,n}) \rightarrow C^{\infty}(\mathrm{Gr}_{i,n}).   \]
\end{enumerate}

Since both the Radon and the cosine transform intertwine the
action of $\mathrm{SO}(n)$, we need some background from harmonic
analysis on Grassmannians to discuss their injectivity
properties. More precisely, we require information on the
decomposition of the space $L^2(\mathrm{Gr}_{i,n})$ into
$\mathrm{SO}(n)$ irreducible subspaces.

Since $\mathrm{SO}(n)$ is a compact Lie group, all its irreducible representations are finite
dimensional. Equivalence classes of irreducible complex
representations of $\mathrm{SO}(n)$ are indexed by their highest
weights (see, e.g., \textbf{\cite[\textnormal{p.\ 219}]{broecker_tomdieck}}), which can be identified with $\lfloor n/2 \rfloor$-tuples of integers
$(\lambda_1,\lambda_2,\ldots,\lambda_{\lfloor n/2 \rfloor})$ such
that
\begin{equation} \label{heiwei}
\left \{\begin{array}{ll} \lambda_1 \geq \lambda_2 \geq \ldots \geq \lambda_{\lfloor n/2 \rfloor} \geq 0 & \quad \mbox{for odd }n, \\
\lambda_1 \geq \lambda_2 \geq \ldots \geq \lambda_{n/2-1} \geq
|\lambda_{n/2}| & \quad \mbox{for even }n.
\end{array} \right .
\end{equation}
We use $\Gamma_\lambda$ to denote any isomorphic
copy of an irreducible representation of $\mathrm{SO}(n)$ with
highest weight
$\lambda=(\lambda_1,\lambda_2,\ldots,\lambda_{\lfloor n/2
\rfloor})$.

\vspace{0.4cm}

\noindent {\bf Examples:}
\begin{enumerate}
\item[(a)] The trivial one-dimensional representation of $\mathrm{SO}(n)$ corresponds to the $\mathrm{SO}(n)$ module $\Gamma_{(0,\ldots,0)}$, while
the standard representation of $\mathrm{SO}(n)$ on $\mathbb{R}^n$ is isomorphic to $\Gamma_{(1,0,\ldots,0)}$.
\item[(b)] The decomposition of $L^2(S^{n-1})$ into an orthogonal sum of $\mathrm{SO}(n)$ \linebreak irreducible subspaces is given by
\begin{equation} \label{decompsn1}
L^2(S^{n-1}) = \bigoplus_{k \geq 0} \Gamma_{(k,0,\ldots,0)}.
\end{equation}
It is well known that here the spaces $\Gamma_{(k,0,\ldots,0)}$ are precisely the spaces of spherical harmonics of degree $k$ in dimension $n$.
\item[(c)] For $1 \leq i \leq n - 1$, the space $L^2(\mathrm{Gr}_{i,n})$ is a sum of orthogonal irreducible representations of
$\mathrm{SO}(n)$ with highest weights $(\lambda_1,\ldots,\lambda_{\lfloor n/2 \rfloor})$ satisfying the following two additional conditions
(see, e.g., \textbf{\cite[\textnormal{Theorem 8.49}]{knapp}}):
\begin{equation} \label{decompgri}
\left \{ \begin{array}{l} \lambda_j = 0\, \mbox{ for all } j > \min\{i, n-i\}, \phantom{wwwwwww} \\
\lambda_1, \ldots, \lambda_{\lfloor n/2 \rfloor}\, \mbox{ are all even.} \end{array} \right .
\end{equation}
\end{enumerate}

\vspace{0.3cm}

Note that since both $S^{n-1}$ and $\mathrm{Gr}_{i,n}$ are Riemannian symmetric spaces the respective decompositions of $L^2(S^{n-1})$ and $L^2(\mathrm{Gr}_{i,n})$
into $\mathrm{SO}(n)$ irreducible subspaces are {\it multiplicity free}. In fact, even more can be said. To this end let $\mathrm{GL}(V)$
denote the general linear group of a vector space $V$.

\vspace{0.4cm}

\noindent {\bf Definition} (see \textbf{\cite[\textnormal{p.\ 14}]{takeuchi}}) \emph{Let $G$ be a compact Lie group and let $H$ be a closed subgroup of $G$.
A representation $\rho: G \rightarrow \mathrm{GL}(V)$ is called \emph{spherical} with respect to $H$ if there exists
a non-zero $v \in V$ which is invariant under $H$, that is, $\rho(\vartheta)v = v$ for every $\vartheta \in H$. Let $V^H$ denote
the subspace of $V$ consisting of all $H$-invariant $v \in V$.}

\vspace{0.4cm}

For the following well-known facts we also refer to \textbf{\cite[\textnormal{p.\ 17}]{takeuchi}}:

\pagebreak

\begin{theorem} \label{thmspher} Let $G$ be a compact Lie group and $H \subseteq G$ a closed subgroup.
\begin{enumerate}
\item[(i)] Every subrepresentation of $L^2(G/H)$ is spherical with respect to $H$.
\item[(ii)] Every irreducible representation which is spherical with respect to $H$ is isomorphic to a subrepresentation of $L^2(G/H)$.
\item[(iii)] If $G/H$ is a Riemannian symmetric space, then $\mathrm{dim}\,V^H = 1$ for every irreducible representation which is spherical with respect to $H$.
\end{enumerate}
\end{theorem}

In order to see how we will use Theorem \ref{thmspher}, consider
the subspace $L^2(\mathrm{Gr}_{i,n})^{\mathrm{SO}(n-1)}$ of
$\mathrm{SO}(n-1)$ invariant functions in
$L^2(\mathrm{Gr}_{i,n})$. By definition, \linebreak every
irreducible subspace of $L^2(\mathrm{Gr}_{i,n})$ which has
non-trivial intersection with
$L^2(\mathrm{Gr}_{i,n})^{\mathrm{SO}(n-1)}$ corresponds to a
spherical subrepresentation with respect to $\mathrm{SO}(n-1)$. Thus, by Theorem
\ref{thmspher} (ii), (\ref{decompsn1}), and (\ref{decompgri}), we conclude that
\begin{equation} \label{corspher}
L^2(\mathrm{Gr}_{i,n})^{\mathrm{SO}(n-1)} \subseteq \bigoplus_{k \geq 0} \Gamma_{(2k,0,\ldots,0)}.
\end{equation}

Let $G$ be a compact Lie group and let $\Gamma$ be a (not necessarily irreducible) finite dimensional complex
$G$ module. Recall that the dual
representation is defined on the dual space $\Gamma^*$ by
\[(\vartheta\,u^*)(v) = u^*(\vartheta^{-1} v), \qquad \vartheta \in G, u^* \in \Gamma^*, v \in \Gamma,  \]
and that $\Gamma$ is called {\it self-dual} if $\Gamma$ and $\Gamma^*$
are isomorphic representations. The module $\Gamma$ is called {\it
real} if there exists a non-degenerate symmetric $G$
invariant bilinear form on $\Gamma$. In particular, if $\Gamma$ is
real, then $\Gamma$ is also self-dual.

Due to (\ref{commute}), the following auxiliary result is crucial to our investigations:

\begin{lem} \label{crit17} Let $H_1$ and $H_2$ be closed subgroups of $\mathrm{SO}(n)$ such that both $\mathrm{SO}(n)/H_1$ and $\mathrm{SO}(n)/H_2$ are Riemannian symmetric
spaces and let
\[L^2(\mathrm{SO}(n)/H_1)=\bigoplus_{\lambda} \Gamma^{(1)}_{\lambda} \qquad \mbox{and} \qquad L^2(\mathrm{SO}(n)/H_2)=\bigoplus_{\xi} \Gamma^{(2)}_{\xi},  \]
where both sums range over suitable equivalence classes of irreducible $\mathrm{SO}(n)$ representations. If $f \in \Gamma^{(1)}_{\lambda}$ is $H_2$ invariant and $\Gamma^{(1)}_{\lambda}$ is real, then
$\widehat{f} \in \Gamma^{(2)}_{\lambda}$ and $\widehat{f}$ is $H_1$ invariant.
\end{lem}
{\it Proof.} Let $L^2(\mathrm{SO}(n))=\bigoplus_{\lambda}\Gamma(\lambda)$ be the decomposition of $L^2(\mathrm{SO}(n))$ into isotypical components (see, e.g., \textbf{\cite[\textnormal{p.\ 70}]{broecker_tomdieck}}). Consider the action of the group $\mathrm{SO}(n) \times \mathrm{SO}(n)$ on $L^2(\mathrm{SO}(n))$ given by
\[((\vartheta,\eta)f)(\zeta) = f(\vartheta^{-1}\zeta\eta).  \]
It is well known (see, e.g., \textbf{\cite[\textnormal{Theorem 1.1}]{takeuchi}}) that
\[\Gamma(\lambda) \cong \Gamma_{\lambda} \otimes \Gamma_{\lambda}^*  \]
as $\mathrm{SO}(n) \times \mathrm{SO}(n)$ modules. If $\Gamma_{\lambda}$ is real, then $\Gamma_{\lambda} \cong \Gamma_{\lambda}^*$ and the isomorphism
$\Theta: \Gamma_{\lambda} \otimes \Gamma_{\lambda} \rightarrow \Gamma(\lambda)$ is given by
\[\Theta(g \otimes h)(\vartheta) = \langle g, \vartheta h \rangle,  \]
where $\langle \,\, ,\, \rangle$ denotes the non-degenerate, symmetric, and $\mathrm{SO}(n)$
invariant \linebreak bilinear form on $\Gamma_{\lambda}$. Since
\[\Theta(g \otimes h)(\vartheta^{-1})=\Theta(h \otimes g)(\vartheta),  \]
it follows that $\widehat{f} \in \Gamma(\lambda)$ for every $f \in
\Gamma(\lambda)$. But, since clearly,
\[\Gamma^{(1)}_{\lambda}, \Gamma^{(2)}_{\lambda} \subseteq \Gamma(\lambda),  \]
we deduce that if $f \in \Gamma^{(1)}_{\lambda}$ is (left) $H_2$ invariant, then $\widehat{f}$ is right $H_2$ invariant and, thus, $\widehat{f} \in \Gamma^{(2)}_{\lambda}$. \hfill $\blacksquare$

\vspace{0.4cm}

We return now to cosine and Radon transforms. Let $1 \leq i \neq j
\leq n - 1$ and let $\Gamma_{\lambda} \subseteq
L^2(\mathrm{Gr}_{i,n})$ be an $\mathrm{SO}(n)$ irreducible
subspace. Since both transforms $C_i: L^2(\mathrm{Gr}_{i,n})
\rightarrow L^2(\mathrm{Gr}_{i,n})$ and $R_{i,j}:
L^2(\mathrm{Gr}_{i,n}) \rightarrow L^2(\mathrm{Gr}_{j,n})$
intertwine the action of $\mathrm{SO}(n)$, the spaces
\[C_i\Gamma_{\lambda}\subseteq L^2(\mathrm{Gr}_{i,n}) \qquad \mbox{and} \qquad R_{i,j}\Gamma_{\lambda} \subseteq L^2(\mathrm{Gr}_{j,n})  \]
are also $\mathrm{SO}(n)$ irreducible subspaces of
$L^2(\mathrm{Gr}_{i,n})$ and $L^2(\mathrm{Gr}_{j,n})$,
respectively. By Schur's Lemma, these subspaces are either
trivial or we must have
\[C_i\Gamma_{\lambda} = \Gamma_{\lambda} \qquad \mbox{and} \qquad R_{i,j}\Gamma_{\lambda} \cong \Gamma_{\lambda}. \]
Therefore the cosine transform $C_i$ must act as a multiple of the
identity on $\Gamma_{\lambda}$, that is, there exist {\it
multipliers} $\mathrm{c}^i_{\lambda} \in \mathbb{R}$ such that
for every $f \in \Gamma_{\lambda}$,
\[C_if = \mathrm{c}^i_{\lambda}\,f. \]
Problems concerning injectivity of $C_i$ are now reduced to
questions as to which of these multipliers are zero.

\pagebreak

By well known classical facts (see, e.g.,
\textbf{\cite{groemer96}}), the {\it spherical} cosine transform
$C_1=C_{n-1}$ and the {\it spherical} Radon transform
$R_{1,n-1}=R_{n-1,1}$ are injective on $L^2(\mathrm{Gr}_1)$ and
$L^2(\mathrm{Gr}_{n-1})$, respectively. In particular,
\begin{equation} \label{1nmin1}
\mathrm{c}^1_{(2k,0,\ldots,0)}\neq 0
\end{equation}
for every $k \in \mathbb{N}$. Moreover, when restricted to smooth
functions these \linebreak spherical transforms are {\it
bijective}.

For general $1 < i < n - 1$, the cosine transform $C_i$ is not
injective as was first shown by Goodey and Howard
\textbf{\cite{goodeyhoward}} (see also, \textbf{\cite{AlBern}}).
However, of particular importance for us is the fact that not
only (\ref{1nmin1}) holds, but that for all $1 \leq i \leq n - 1$
and for every $k \in \mathbb{N}$,
\begin{equation} \label{cimult}
\mathrm{c}^i_{(2k,0,\ldots,0)}=\frac{n \kappa_i \kappa_{n-i}}{2 \kappa_{n-1}}
 {n \choose i}^{-1}\mathrm{c}^1_{(2k,0,\ldots,0)} \neq 0.
\end{equation}
This relation was first obtained by Goodey and Zhang
\textbf{\cite[\textnormal{Lemma 2.1}]{goodeyzhang}}. \linebreak
Using (\ref{cimult}) one can show that the restriction of $C_i$
to the subspace of $C^{\infty}(\mathrm{Gr}_{i,n})$ defined by
\begin{equation} \label{subspace}
\mathrm{cl}_{C^{\infty}} \bigoplus_{k \geq 0} \Gamma_{(2k,0,\ldots,0)}
\end{equation}
is {\it bijective}. Here $\mathrm{cl}_{C^{\infty}}$ denotes the
closure in the $C^{\infty}$ topology.

If $1 \leq i < j \leq n - 1$, then $R_{i,j}$ is injective if and
only if $i + j \leq n$, whereas if $i > j$, then $R_{i,j}$ is
injective if and only if $i + j \geq n$, see Grinberg
\textbf{\cite{grinberg86}}. Moreover, it was shown in
\textbf{\cite{goodeyzhang}} that for all $1 \leq i \neq j \leq n -
1$ also the restriction of the Radon transform $R_{i,j}$ to the
subspace defined by (\ref{subspace}) is bijective.

The following consequences of these facts and Lemma \ref{crit17}
will be the key ingredients in the proof of Theorems 1 and 2.

\begin{lem} \label{keylem} Suppose that $1 \leq i \leq n - 1$.
\begin{enumerate}
\item[(i)] If $\mu$ is an $\mathrm{S}(\mathrm{O}(i) \times \mathrm{O}(n-i))$ invariant
measure on $S^{n-1}$, then
\[\left ( \widehat{C_i \widehat{\mu}}\, \right )\!(u) = \frac{n \kappa_i \kappa_{n-i}}{2 \kappa_{n-1}}
 {n \choose i}^{-1}\!\! \int_{S^{n-1}} | u \cdot v|\,d\mu(v), \qquad u \in S^{n-1}.  \]
\item[(ii)] If $f \in C^{\infty}(S^{n-1})$ is $\mathrm{S}(\mathrm{O}(i) \times
\mathrm{O}(n-i))$ invariant, then there exists a unique zonal $g
\in C^{\infty}(S^{n-1})$ such that
\begin{equation} \label{keylem2}
f = \lambda_{i,n-1} \ast g.
\end{equation}
\end{enumerate}
\end{lem}

\noindent {\it Proof.} In order to prove (i) we may assume that $d\mu= h\, du$ for some
$\mathrm{S}(\mathrm{O}(i) \times \mathrm{O}(n-i))$ invariant function $h \in C(S^{n-1})$.
Let $h = \sum_{k \geq 0} h_k$ be the
decomposition of $h$ into spherical harmonics, that is, $h_k \in
\Gamma_{(k,0,\ldots,0)}$. Since $h$ is $\mathrm{S}(\mathrm{O}(i)
\times \mathrm{O}(n-i))$ invariant, so is each $h_k$.

It is well known (cf.\
\textbf{\cite[\textnormal{p.292}]{broecker_tomdieck}}) that each
$\Gamma_{(k,0,\ldots,0)}$ is a real $\mathrm{SO}(n)$ module.
Thus, since $\widehat{h} \in C(\mathrm{Gr}_{i,n})^{\mathrm{SO}(n-1)}$, it follows from  Lemma \ref{crit17} and  \eqref{decompgri} that on one hand $h_{2k+1}=0$ and on the other hand $\widehat{h}_{2k} \in \Gamma_{(2k,0,\ldots,0)}$ for every $k \in
\mathbb{N}$. Hence, by (\ref{cimult}) we obtain
\[\widehat{C_i \widehat{h}} = \sum_{k\geq 0} \widehat{C_i \widehat{h}_{2k}} = \frac{n \kappa_i \kappa_{n-i}}{2 \kappa_{n-1}}
 {n \choose i}^{-1} \sum_{k \geq 0} C_1 h_{2k}= \frac{n \kappa_i \kappa_{n-i}}{2 \kappa_{n-1}}
 {n \choose i}^{-1} C_1h, \]
which is precisely the claim from (i).

\vspace{0.1cm}

In order to prove (ii), first note that for any zonal $g \in C^{\infty}(S^{n-1})$ we have $g = \widehat{g}$ (cf. \textbf{\cite{Schu06a}}). Thus,
by (\ref{commute}) and (\ref{radconv}), relation (\ref{keylem2}) is equivalent to
\[\widehat{f} = g \ast \widehat{\lambda}_{i,n-1}=R_{n-1,i}\,g.  \]
Since $\widehat{f} \in C^{\infty}(\mathrm{Gr}_{i,n})^{\mathrm{SO}(n-1)}$, it follows from (\ref{corspher}) that $\widehat{f}$ is
contained in the subspace of $C^{\infty}(\mathrm{Gr}_{i,n})$ defined in (\ref{subspace}). Since the Radon transform $R_{n-1,i}$ is an $\mathrm{SO}(n)$
intertwining bijection from even smooth functions
on $S^{n-1}$ to this subspace, statement (ii) follows. \hfill $\blacksquare$

\vspace{1cm}

\centerline{\large{\bf{ \setcounter{abschnitt}{4}
\arabic{abschnitt}. Even translation invariant valuations}}}

\reseteqn \alpheqn \setcounter{theorem}{0}

\vspace{0.6cm}

In this section we collect basic facts about convex bodies (see,
e.g., \textbf{\cite{schneider93}}) and the background material
from the theory of {\it even} translation invariant scalar and
convex body valued valuations. In particular, we recall Klain's
\linebreak embedding of even translation invariant valuations into
functions on the Grassmannian, the definitions of Alesker's Hard
Lefschetz operators and his Fourier type transform on translation
invariant valuations.

The definition of the support function $h(K,u)=\max\{u\cdot x: x
\in K\}$, $u \in S^{n-1}$, of a convex body $K \in \mathcal{K}^n$,
implies that $h(\vartheta K,u)=h(K, \vartheta^{-1}u)$ for every $u
\in S^{n-1}$ and $\vartheta \in \mathrm{O}(n)$. For $K_1, K_2 \in
\mathcal{K}^n$ and $\lambda_1, \lambda_2 \geq 0$, the support
function of the Minkowski linear combination
$\lambda_1K_1+\lambda_2K_2$ is given by
\[h(\lambda_1K_1+\lambda_2K_2,\cdot)=\lambda_1h(K_1,\cdot)+\lambda_2h(K_2,\cdot).\]

The {\it surface area measure} $S_{n-1}(K,\cdot)$ of a convex
body $K$ is defined for Borel sets $\omega \subseteq S^{n-1}$, as
the $(n-1)$-dimensional Hausdorff measure of the set of all
boundary points of $K$ at which there exists a normal vector of
$K$ belonging to $\omega$. If the body $K \in \mathcal{K}^n$ has
non-empty interior, then $K$ is determined up to translations by
its surface area measure.

Let $B$ denote the Euclidean unit ball in $\mathbb{R}^n$. The
surface area measure of $K \in \mathcal{K}^n$ satisfies the
Steiner-type formula
\begin{equation} \label{steiner}
S_{n-1}(K + \varepsilon B,\cdot)=\sum \limits_{i=0}^{n-1}
\varepsilon^{n-1-i} {n - 1 \choose i}S_i(K,\cdot).
\end{equation}
The measure $S_i(K,\cdot)$, $0 \leq i \leq n - 1$, is called the
{\it area measure of order $i$} of $K \in \mathcal{K}^n$. The
relation $S_i(\lambda K,\cdot)=\lambda^iS_i(K,\cdot)$ holds for
all $K \in \mathcal{K}^n$ and every $\lambda > 0$. For
$\vartheta \in \mathrm{O}(n)$, we have $S_i(\vartheta
K,\cdot)=\vartheta S_i(K,\cdot)$.

The area measure
$S_i(K,\cdot)$ is (up to normalization) a localization of the
$i$th {\it intrinsic volume} $V_i(K)$ of a convex body $K$. More
precisely,
\[S_i(K,S^{n-1}) = n {n \choose i}^{-1}\kappa_{n-i}V_i(K).   \]

For $i \in \{1, \ldots, n - 1\}$, a special case of the
Cauchy--Kubota formulas implies that
\begin{equation} \label{radsurfcos}
R_{i,n-1}\mathrm{vol}_i(K|\,\cdot\,) =
\frac{\kappa_i}{4\kappa_{n-1}}\, C_{n-1}\left(S_i(K,\cdot) +
S_i(-K,\cdot) \right ),
\end{equation}
where we identify the even measure $S_i(K,\cdot) + S_i(-K,\cdot)$
on $S^{n-1}$ with a measure on $\mathrm{Gr}_{n-1,n}$ of the same total mass.

Intrinsic volumes and area measures are both valuations. In the case of
intrinsic volumes we have scalar-valued
valuations; in the case of area measures we have valuations with values in the set of Borel
measures on \nolinebreak $S^{n-1}$.

If $G$ is a group of affine transformations on $\mathbb{R}^n$, a
valuation $\phi$ is called \linebreak $G$-invariant if $\phi(gK) =
\phi(K)$ for all $K \in \mathcal{K}^n$ and every $g \in G$. We
denote the vector space of continuous translation invariant
scalar-valued valuations by $\mathbf{Val}$. A seminal result in
the structure theory of translation invariant valuations was
obtained by McMullen \textbf{\cite{McMullen77}}, who showed that
\begin{equation} \label{mcmullen}
\mathbf{Val} = \bigoplus_{0\leq i \leq n} \mathbf{Val}_i  =  \bigoplus_{0\leq i \leq n} \left ( \mathbf{Val}_i^+ \oplus
\mathbf{Val}_i^-\right ),
\end{equation}
where $\mathbf{Val}_i^+ \subseteq \mathbf{Val}$ denotes the
subspace of {\it even} valuations (homogeneous) of degree $i$, $\mathbf{Val}_i^-$ denotes the subspace of {\it odd} valuations of
degree $i$, and $\mathbf{Val}_i=\mathbf{Val}_i^+ \oplus
\mathbf{Val}_i^-$ is the subspace of valuations of degree $i$.

\pagebreak

From (\ref{mcmullen}), it follows easily that the
space $\mathbf{Val}$ becomes a Banach space, when endowed with
the norm $\|\phi \|=\sup \{|\phi(K)|: K \subseteq B\}$.
The general linear group $\mathrm{GL}(n)$ acts on this Banach
space in a natural way: For every $A \in
\mathrm{GL}(n)$ and every $K \in \mathcal{K}^n$,
\[(A\phi)(K)=\phi(A^{-1}K), \qquad \phi \in \mathbf{Val}.  \]

Assume that $1 \leq i \leq n - 1$ and define for any finite Borel
measure $\mu$ on $\mathrm{Gr}_{i,n}$ an even valuation
$\mathrm{Cr}_i\mu \in \mathbf{Val}_i^+$ by
\[(\mathrm{Cr}_i\mu)(K) = \int_{\mathrm{Gr}_{i,n}} \mathrm{vol}_i(K|E)\,d\mu(E).   \]
It follows from a deep result of Alesker
\textbf{\cite{Alesker01}} that the image of the map
$\mathrm{Cr}_i$ is dense in $\mathbf{Val}_i^+$.
 This leads to the following

\vspace{0.3cm}

\noindent {\bf Definition} \emph{A finite Borel measure $\mu$ on
$\mathrm{Gr}_{i,n}$, $1 \leq i \leq n - 1$, is called a
\emph{Crofton measure} for the valuation $\phi \in
\mathbf{Val}_i^+$ if $\mathrm{Cr}_i\mu = \phi$.}

\vspace{0.3cm}

Using the notion of smooth vectors of a representation (see \textbf{\cite[\textnormal{p.\ 31}]{wallach1}}) and an
imbedding of the space $\mathbf{Val}_i^+$ in
$C(\mathrm{Gr}_{i,n})$ by Klain, a more precise description of
valuations admitting a Crofton measure is possible.

\vspace{0.3cm}

\noindent {\bf Definition} \emph{A valuation $\phi \in
\mathbf{Val}$ is called smooth if the map $\mathrm{GL}(n)
\rightarrow \mathbf{Val}$ defined by $A \mapsto A\phi$ is
infinitely differentiable.}

\vspace{0.3cm}

We write $\mathbf{Val}^{\infty}$ for the space of smooth
translation invariant valuations, and we use
$\mathbf{Val}_i^{\infty}$ and $\mathbf{Val}_i^{\pm,\infty}$ for the subspaces of smooth
valuations in $\mathbf{Val}_i$ and $\mathbf{Val}_i^{\pm}$, respectively. It is well known (cf.\
\textbf{\cite[\textnormal{p.\ 32}]{wallach1}}) that the set of
smooth valuations $\mathbf{Val}_i^{\pm,\infty}$ is a dense
$\mathrm{GL}(n)$ invariant subspace of $\mathbf{Val}_i^{\pm}$.

For $1 \leq i \leq n - 1$, the {\it Klain map} $\mathrm{Kl}_i:
\mathbf{Val}_i^+ \rightarrow C(\mathrm{Gr}_{i,n})$, $\phi \mapsto \mathrm{Kl}_i\phi$, is defined as
follows: For $\phi \in \mathbf{Val}_i^+$ and every $E \in
\mathrm{Gr}_{i,n}$, consider the restriction $\phi_E$ of $\phi$
to convex bodies in $E$. This is a continuous translation
invariant valuation of degree $i$ in $E$. Thus, by a result of
Hadwiger \textbf{\cite[\textnormal{p.\ 79}]{hadwiger51}}, \linebreak $\phi_E
= (\mathrm{Kl}_i\phi)(E)\,\mathrm{vol}_i$, where
$(\mathrm{Kl}_i\phi)(E)$ is a constant depending only on $E$. This gives
rise to a continuous function $\mathrm{Kl}_i\phi \in C(\mathrm{Gr}_{i,n})$, called the \emph{Klain
function} of the valuation $\phi$. By an important
result of Klain \textbf{\cite{klain00}}, the Klain map $\mathrm{Kl}_i$ is injective for every $i \in \{1, \ldots, n - 1\}$.

Consider now the restriction of the map $\mathrm{Cr}_i$, $1 \leq i
\leq n - 1$, to smooth functions:
\begin{equation*}
(\mathrm{Cr}_if)(K)=\int_{\mathrm{Gr}_{i,n}}\mathrm{vol}_i(K|E)f(E)\,dE,
\qquad f \in C^{\infty}(\mathrm{Gr}_{i,n}).
\end{equation*}
It is not difficult to see that the valuation $\mathrm{Cr}_if$ is
smooth, i.e., $\mathrm{Cr}_if \in \mathbf{Val}_i^{+,\infty}$.

\noindent Moreover, if $F \in \mathrm{Gr}_{i,n}$, then, for any
$f \in C^{\infty}(\mathrm{Gr}_{i,n})$ and convex body $K \subseteq
F$,
\begin{equation*}
(\mathrm{Cr}_if)(K)=\mathrm{vol}_i(K)\int_{\mathrm{Gr}_{i,n}}
|\cos(E,F)|f(E)\,dE.
\end{equation*}
Thus, the Klain function of the valuation $\mathrm{Cr}_if$ is
equal to the cosine transform $C_if$ of $f$. From the
main result of \textbf{\cite{AlBern}}, Alesker
\textbf{\cite[\textnormal{p.\ 73}]{Alesker03}} deduced the
following:

\begin{theorem} \label{klaincoscrof} \emph{(Alesker and Bernstein \textbf{\cite{AlBern}}, Alesker \textbf{\cite{Alesker03}})}  Let $1 \leq i \leq n - 1$. The image of the
Klain map $\mathrm{Kl}_i: \mathbf{Val}_i^{+,\infty} \rightarrow
C^{\infty}(\mathrm{Gr}_{i,n})$ coincides with the image of the
cosine transform $C_i: C^{\infty}(\mathrm{Gr}_{i,n}) \rightarrow
C^{\infty}(\mathrm{Gr}_{i,n})$. Moreover, for every valuation $\phi
\in \mathbf{Val}_i^{+,\infty}$, there exists a smooth Crofton measure.
\end{theorem}

We turn now to Alesker's Hard Lefschetz operators. It is well known that
McMullen's decomposition (\ref{mcmullen}) of $\mathbf{Val}$ into subspaces of homogeneous valuations
implies that for every $\phi \in \mathbf{Val}$ and all $K \in \mathcal{K}^n$ the Steiner type formula
\begin{equation} \label{steinertype1}
\phi(K + rB)=\sum_{j=0}^n r^{n-j} \phi^{(j)}(K),
\end{equation}
where $\phi^{(j)} \in \mathbf{Val}$ for $0 \leq j \leq n$, holds
for every $r \geq 0$. Note that $\phi^{(j)}$ is in general not homogeneous.

In turn, (\ref{steinertype1}) gives rise to a derivation operator $\Lambda: \mathbf{Val} \rightarrow \mathbf{Val}$ defined by
\[(\Lambda \phi)(K) = \left . \frac{d}{dt} \right |_{t=0} \phi(K + tB).  \]
Note that $\Lambda$ commutes with the action of $\mathrm{O}(n)$ and that it preserves parity.
Moreover, if $\phi \in \mathbf{Val}_i$, then $\Lambda \phi \in \mathbf{Val}_{i-1}$.

The importance of the operator $\Lambda$ becomes evident from
the following Hard Lefschetz type theorem established for even
valuations by Alesker \textbf{\cite{Alesker03}} and for general
valuations by Bernig and Br\"ocker \textbf{\cite{Bernig07b}}:

\begin{theorem} \label{hardlef1}  \emph{(Alesker \textbf{\cite{Alesker03}}, Bernig and Br\"ocker \textbf{\cite{Bernig07b}})} Let $1 \leq i \leq n$.
\begin{enumerate}
\item[(i)] The operator $\Lambda: \mathbf{Val}^{\infty}_{i} \rightarrow
\mathbf{Val}^{\infty}_{i-1}$ is injective if $2i-1 \geq n$ and
surjective if $2i - 1 \leq n$.
\item[(ii)] If $2i \geq n$, then $\Lambda^{2i-n}:
\mathbf{Val}^{\infty}_{i} \rightarrow
\mathbf{Val}^{\infty}_{n-i}$ is an isomorphism.
\end{enumerate}
\end{theorem}

More recently, a dual version of Theorem \ref{hardlef1} was established by Alesker, in \textbf{\cite{Alesker04}} for even valuations and in \textbf{\cite{alesker10}} for general ones. There, the
derivation operator $\Lambda$ is replaced by an integral operator $\mathfrak{L}: \mathbf{Val} \rightarrow \mathbf{Val}$
defined by
\begin{equation} \label{hardlefint}
(\mathfrak{L}\phi)(K) = \int_{\mathrm{AGr}_{n-1,n}} \phi(K \cap E)\,d\sigma_{n-1}(E).
\end{equation}
Here and in the following $\mathrm{AGr}_{i,n}$, denotes the affine Grassmannian of $i$ planes in $\mathbb{R}^n$ and $\sigma_{i}$
is the invariant measure on $\mathrm{AGr}_{i,n}$ normalized such that the set of planes having non-empty intersection with the Euclidean unit ball has measure
\[\left[\begin{array}{c} n\\i \end{array}\right]\kappa_{n-i} := {n \choose i} \frac{\kappa_n}{\kappa_i}.  \]

The operator $\mathfrak{L}$ was originally introduced by Alesker in a different way in connection
with a newly discovered product structure on the space $\mathbf{Val}$. Only in \textbf{\cite{Bernig07}} Bernig
showed that the original definition coincides with (\ref{hardlefint}). We also note that $\mathfrak{L}$ commutes with the action of $\mathrm{O}(n)$ and that it preserves parity.
Moreover, if $\phi \in \mathbf{Val}_i$, then $\mathfrak{L}\phi \in \mathbf{Val}_{i+1}$.

The Hard Lefschetz type theorem for the operator $\mathfrak{L}$ is

\begin{theorem} \label{hardlef2}  \emph{(Alesker \textbf{\cite{alesker10}})} Let $0 \leq i \leq n - 1$.
\begin{enumerate}
\item[(i)] The operator $\mathfrak{L}: \mathbf{Val}^{\infty}_{i} \rightarrow
\mathbf{Val}^{\infty}_{i+1}$ is injective if $2i+1 \leq n$ and
surjective if $2i + 1 \geq n$.
\item[(ii)] If $2i \leq n$, then $\mathfrak{L}^{n-2i}:
\mathbf{Val}^{\infty}_{i} \rightarrow
\mathbf{Val}^{\infty}_{n-i}$ is an isomorphism.
\end{enumerate}
\end{theorem}

The dual nature of Theorems \ref{hardlef1} and \ref{hardlef2} was the basis not only for the proof of
Theorem \ref{hardlef2} but for the discovery of another fundamental duality transform, now called
the Alesker--Fourier transform (which shares various formal similarities with the classical Fourier transform). In fact,
it was shown by Bernig and Fu \textbf{\cite{bernigfu06}} for even valuations and recently by Alesker \textbf{\cite{alesker10}} for general valuations
that both versions of the Hard Lefschetz Theorem
are equivalent via the Alesker--Fourier transform.

We will define the Alesker--Fourier transform $\mathbb{F}: \mathbf{Val}^{\infty} \rightarrow \mathbf{Val}^{\infty}$ only for even valuations
and refer to the article \textbf{\cite{alesker10}} for the odd case which is much more involved and will not be needed in the following: If
$\phi \in \mathbf{Val}^{+,\infty}_i$, $1 \leq i \leq n - 1$, then $\mathbb{F}\phi \in \mathbf{Val}^{+,\infty}_{n-i}$ is the valuation whose
Klain function is given by
\begin{equation} \label{deffourier}
\mathrm{Kl}_{n-i}(\mathbb{F}\phi) = (\mathrm{Kl}_i\phi)^{\bot}.
\end{equation}
In order to see that the linear operator $\mathbb{F}$ is well defined, use (\ref{cifbot}) and Theorem \ref{klaincoscrof}. Clearly, $\mathbb{F}$ is an involution that commutes with the action
of $\mathrm{O}(n)$.
The derivation operator $\Lambda$ and the integral operator $\mathfrak{L}$ are related by
\begin{equation} \label{hardfourier}
\mathbb{F} \circ \Lambda = 2\, \mathfrak{L} \circ \mathbb{F}.
\end{equation}
This was first observed by Bernig and Fu \textbf{\cite{bernigfu06}} for even valuations and proved in general by Alesker in \textbf{\cite{alesker10}}.

\vspace{0.2cm}

We conclude this section by collecting previously obtained
results on translation invariant and $\mathrm{SO}(n)$ equivariant Minkowski valuations needed in the next section.
To this end we denote by $\mathbf{MVal}^{(+)}$ the set of continuous translation invariant (even) Minkowski valuations, and
we write $\mathbf{MVal}_i^{(+)}$, $0 \leq i \leq n$, for its subset of all (even) Minkowski valuations of degree $i$.

A Minkowski valuation $\Phi \in \mathbf{MVal}$ is called \emph{$\mathrm{SO}(n)$ equivariant} if for all $K \in \mathcal{K}^n$
and every $\vartheta \in \mathrm{SO}(n)$,
\[\Phi(\vartheta K) = \vartheta \Phi K.  \]
For a number of well-known examples of Minkowski valuations that are $\mathrm{SO}(n)$ equivariant we refer to the
articles \textbf{\cite{kiderlen05, Ludwig:Minkowski, Schu09}} and the next section.

Recently, Parapatits and the second author \textbf{\cite{parapwann}} have shown that, in \linebreak general, a decomposition of a Minkowski valuation $\Phi \in \mathbf{MVal}$
into a sum of homogeneous Minkowski valuations is not possible (cf.\ also \textbf{\cite{parapschu}}). However, from an application of McMullen's
decomposition (\ref{mcmullen}) it is possible to deduce the following: If $\Phi \in \mathbf{MVal}^{(+)}$, then there exist convex bodies $L_0, L_n
 \in \mathcal{K}^n$ and for every $K \in \mathcal{K}^n$, (even) functions $g_i(K,\cdot) \in C(S^{n-1})$ such that
\begin{equation} \label{suppdecomp}
h(\Phi K,\cdot) = h(L_0,\cdot) + \sum_{i=1}^{n-1} g_i(K,\cdot) +
V(K)h(L_n,\cdot).
\end{equation}
Moreover, for each $i \in \{1, \ldots, n - 1\}$:
\begin{enumerate}
\item[(i)] The map $K \mapsto g_i(K,\cdot)$ is a continuous translation
invariant valuation of degree $i$.
\item[(ii)] If $\Phi$ is $\mathrm{SO}(n)$ equivariant, then $L_0$ and $L_n$ are Euclidean balls and for every $\vartheta \in \mathrm{SO}(n)$ and $K \in \mathcal{K}^n$ we have
$g_i(\vartheta K,u) = g_i(K,\vartheta^{-1}u)$.
\end{enumerate}
In general, the functions $g_i(K,\cdot)$ need not be support functions (see \textbf{\cite{parapwann}}). However, it is an important open problem
whether for $\mathrm{SO}(n)$ equivariant $\Phi \in \mathbf{MVal}$, the $g_i(K,\cdot)$ are support functions for every $i \in \{1, \ldots, n - 1\}$.

\pagebreak

If $\Phi \in \mathbf{MVal}$ is $\mathrm{SO}(n)$ equivariant, then for $\bar{\vartheta} \in S^{n-1}$ we have
\[h(\Phi K,\bar{\vartheta}) = h(\Phi K,\vartheta \bar{e})=h(\vartheta^{-1}(\Phi K),\bar{e})=h(\Phi(\vartheta^{-1}K),\bar{e}).   \]
Consequently, the real valued valuation $K \mapsto h(\Phi K,\bar{e})$ uniquely determines the Minkowski valuation $\Phi$.
This motivates the following:

\vspace{0.3cm}

\noindent {\bf Definition} \emph{Suppose that $\Phi \in \mathbf{MVal^{(+)}}$ is $\mathrm{SO}(n)$ equivariant. The $\mathrm{SO}(n - 1)$
invariant real valued valuation $\varphi \in \mathbf{Val}^{(+)}$, defined by}
\[\varphi(K) = h(\Phi K,\bar{e}), \qquad K \in \mathcal{K}^n,  \]
\emph{is called the \emph{associated real valued valuation} of $\Phi \in \mathbf{MVal^{(+)}}$. We say that $\Phi$ is \emph{smooth} if its  associated real valued valuation $\varphi$ is smooth. }

\vspace{0.3cm}

The notion of smoothness for translation invariant and $\mathrm{SO}(n)$ equivariant Minkowski valuations was introduced by the first author in \textbf{\cite{Schu09}}.
There, it was also shown that any \emph{even} $\Phi \in \mathbf{MVal}$ which is $\mathrm{SO}(n)$ equivariant can
be approximated uniformly on compact subsets of $\mathcal{K}^n$ by smooth ones.

The following is a reformulation (and slight variation) of Theorem 1 stated in Section 2:

\begin{theorem} \label{mainduke} {\bf(\!\!\textbf{\cite{Schu09}})} Suppose that $i \in \{1, \ldots, n - 1\}$. If $\Phi_i \in \mathbf{MVal}_i^{+}$ is $\mathrm{SO}(n)$ equivariant and smooth,
then there exists an $\mathrm{S}(\mathrm{O}(i) \times \mathrm{O}(n-i))$ invariant $f \in C^{\infty}(S^{n-1})$ such that for every $K \in \mathcal{K}^n$,
\[h(\Phi_iK,\cdot) = \mathrm{vol}_i(K|\cdot) \ast f.   \]
\end{theorem}

Theorem \ref{mainduke} was proved in \textbf{\cite{Schu09}} for Minkowski valuations which are \linebreak $\mathrm{O}(n)$ equivariant. However,
in \textbf{\cite[\textnormal{Lemma 7.1}]{ABS2011}} it was shown that any \linebreak $\mathrm{SO}(n)$ equivariant Minkowski valuation is
also $\mathrm{O}(n)$ equivariant. Another difference between Theorem \ref{mainduke} and the corresponding statement in
\textbf{\cite{Schu09}} is that the convolution appearing in Theorem \ref{mainduke} is induced from $\mathrm{SO}(n)$ while in \textbf{\cite{Schu09}}
the convolution used was induced from $\mathrm{O}(n)$. However, since the convolution of functions on $\mathrm{Gr}_{i,n}$ with $\mathrm{O}(i) \times \mathrm{O}(n-i)$ invariant functions on $S^{n-1}$ does not depend on whether it is induced from $\mathrm{SO}(n)$ or $\mathrm{O}(n)$, the result from \textbf{\cite{Schu09}} implies Theorem \ref{mainduke}.
Conversely, since every $\mathrm{S}(\mathrm{O}(i) \times \mathrm{O}(n-i))$ invariant function on $S^{n-1}$ is
also $\mathrm{O}(i) \times \mathrm{O}(n-i)$ invariant (cf.\ the Appendix), the result from \textbf{\cite{Schu09}} follows from Theorem \ref{mainduke}.
In particular, Theorem 2 also implies the uniqueness of (spherical) Crofton measures in the sense of \textbf{\cite{Schu09}}.

\vspace{0.1cm}

A final result from \textbf{\cite{Schu09}} which we will need concerns the Klain function of the real valued valuation associated
with a translation invariant and $\mathrm{SO}(n)$ equivariant even Minkowski valuation:

\begin{theorem} \label{thmklainbod} {\bf(\!\!\textbf{\cite{Schu09}})} Suppose that $i \in \{1, \ldots, n - 1\}$. If $\Phi_i \in \mathbf{MVal}_i^{+}$ is $\mathrm{SO}(n)$ equivariant and
$\varphi_i \in \mathbf{Val}^+_i$ denotes its associated real valued valuation, then there exists a unique $\mathrm{O}(i) \times \mathrm{O}(n-i)$ invariant
convex body $M \in \mathcal{K}^n$, called the \emph{Klain body} of $\Phi_i$, such that
\[\widehat{\mathrm{Kl}_i\varphi_i} = h(M,\cdot).   \]
\end{theorem}

We emphasize that the Klain body of $\Phi_i$ determines the valuation $\Phi_i$ uniquely.
It is easy to see (and was proved in \textbf{\cite{Schu09}}) that
\[M = \Phi_iK_{\bar{e}},   \]
where $K_{\bar{e}}$ is any convex body in $\bar{e} \in \mathrm{Gr}_{i,n}$ such that $\mathrm{vol}_i(K_{\bar{e}}) = 1$.

\vspace{1cm}

\centerline{\large{\bf{ \setcounter{abschnitt}{5}
\arabic{abschnitt}. Proof of the main results}}}

\reseteqn \alpheqn \setcounter{theorem}{0}

\vspace{0.6cm}

After these preparations we are now in a position to prove our
main \linebreak results, Theorem 2 and Theorem 3 from Section 2. At the end of the section we also
provide explicit integral transforms which relate the previously known
representations of a translation invariant and
$\mathrm{SO}(n)$ equivariant even Minkowski valuation with our new one.

We begin with the proof of Theorem 2 from Section 2:

\vspace{0.3cm}

\noindent {\it Proof of Theorem 2.} Let $\varphi_i \in \mathbf{Val}_i^+$ denote the associated real valued valuation of $\Phi_i$. From
\[h(\Phi_i K,\cdot) = \mathrm{vol}_i(K|\cdot) \ast \mu,  \]
it follows easily (cf.\ \textbf{\cite[\textnormal{p.\ 19}]{Schu09}}) that
\begin{equation} \label{proof1}
\varphi_i(K)=\int_{\mathrm{Gr}_{i,n}} \mathrm{vol}_i(K|E)\,d\widehat{\mu}(E)
\end{equation}
for every $K \in \mathcal{K}^n$. Note that the measure $\widehat{\mu}$ on $\mathrm{Gr}_{i,n}$ is $\mathrm{SO}(n - 1)$ invariant.

From (\ref{proof1}), the remarks before Theorem \ref{klaincoscrof}, and Theorem \ref{thmklainbod}, we obtain
\begin{equation} \label{proof2}
\widehat{\mathrm{Kl}_i\varphi_i} = \widehat{C_i\widehat{\mu}} =
h(M,\cdot).
\end{equation}

\pagebreak

\noindent Here the $\mathrm{O}(i) \times \mathrm{O}(n-i)$ invariant
convex body $M \in \mathcal{K}^n$ is the Klain body of $\Phi_i$ which is uniquely determined by $\Phi_i$. If we define
\[L = \frac{2\kappa_{n-1}}{n\kappa_i\kappa_{n-i}}{n \choose i}\,M,  \]
then it follows from (\ref{proof2}) and Lemma \ref{keylem} (i) that
\[h(L,u) = \int_{S^{n-1}}|u\cdot v|\,d\mu(v), \qquad u \in S^{n-1}.  \]
Thus, by the injectivity of the spherical cosine transform, it follows that the spherical Crofton measure $\mu$ is uniquely determined by $\Phi_i$.
\hfill $\blacksquare$

\vspace{0.4cm}

By Theorem \ref{main1}, the spherical cosine transform of the spherical Crofton measure of an even Minkowski valuation is
necessarily a support function. It is an important open problem whether this condition is also sufficient for an $\mathrm{O}(i) \times \mathrm{O}(n-i)$
invariant measure on $S^{n-1}$ to be the spherical Crofton measure of an even Minkowski valuation of degree $i$. (If $i = n - 1$, then this is the case, see \textbf{\cite{Schu06a}}.)

\vspace{0.2cm}

Using the uniqueness of spherical Crofton measures, we can now deduce Theorem 3 from Theorem \ref{mainduke} and Lemma \ref{keylem} (ii):

\vspace{0.3cm}

\noindent {\it Proof of Theorem 3.}  First note that $C([-1,1])$ is in one to one
correspondence with the subspace of zonal functions in
$C(S^{n-1})$ via the map $g \mapsto g(\bar{e}\cdot\,.\,)$.
Therefore, for any even zonal function $\breve{g} \in
C^{\infty}(S^{n-1})$, there exists an even function $g \in
C^{\infty}((-1,1)) \cap C([-1,1])$ such that $\breve{g} =
g(\bar{e}\cdot\,.\,)$. Moreover, for every measure $\tau$ on $S^{n-1}$ and every $u \in S^{n-1}$,
\[\int_{S^{n-1}} g(u\cdot v)\,d\tau(v) = (\tau \ast \breve{g})(u).  \]
Thus, in order to prove Theorem \ref{main2}, it is sufficient to show that there exists a unique even zonal function $\breve{g} \in C^{\infty}(S^{n-1})$ such that for every $K \in \mathcal{K}^n$,
\begin{equation} \label{genconv17}
h(\Phi_iK,\cdot) = S_i(K,\cdot) \ast \breve{g}.
\end{equation}

By Theorem \ref{mainduke}, there exists an $\mathrm{S}(\mathrm{O}(i) \times \mathrm{O}(n-i))$ invariant function $f \in C^{\infty}(S^{n-1})$ such that for every $K \in \mathcal{K}^n$,
\begin{equation} \label{proof3}
h(\Phi_iK,\cdot) = \mathrm{vol}_i(K|\,\cdot\,) \ast f.
\end{equation}
Since any $\mathrm{S}(\mathrm{O}(i) \times \mathrm{O}(n-i))$ invariant function is also $\mathrm{O}(i) \times \mathrm{O}(n-i)$ invariant (cf.\ Appendix), we can view $f$ as the (smooth) density of a spherical Crofton measure for
$\Phi_i$. By Theorem \ref{main1}, $f$ is uniquely determined by $\Phi_i$.

From Lemma \ref{keylem} (ii), it follows that there exists a unique even zonal \linebreak $\tilde{g} \in C^{\infty}(S^{n-1})$ such that
\begin{equation} \label{proof4}
f=\lambda_{i,n-1} \ast \tilde{g}.
\end{equation}
Hence, plugging (\ref{proof4}) into (\ref{proof3}) and using (\ref{radconv}), we obtain
\[h(\Phi_iK,\cdot) = \mathrm{vol}_i(K|\cdot) \ast \lambda_{i,n-1} \ast \tilde{g} = R_{i,n-1}\mathrm{vol}_i(K|\cdot) \ast \tilde{g}.  \]
Thus, if we define
\[\breve{g}(u) = \frac{\kappa_i}{2\kappa_{n-1}}\,C_{n-1}\,\tilde{g},  \]
then, using (\ref{radsurfcos}), (\ref{cosconv}) and the fact that the spherical convolution of zonal functions is commutative, we arrive at
\[h(\Phi_iK,\cdot)= \frac{\kappa_i}{2\kappa_{n-1}} S_i(K,\cdot) \ast |\cos(\bar{e},\cdot\,)| \ast \tilde{g} = S_i(K,\cdot) \ast \breve{g}. \]

\vspace{-0.5cm}

\hfill $\blacksquare$

\vspace{0.3cm}

Theorem \ref{main2} gives rise to the following:

\vspace{0.3cm}

\noindent {\bf Definition} \emph{Let $i \in \{1, \ldots, n - 1\}$. We call a zonal function $\breve{g} \in C(S^{n-1})$ (or its associated function $g \in C([-1,1])$)
a \emph{generating function} for $\Phi_i \in \mathbf{MVal}_i$ if (\ref{genconv17}) holds for every $K \in \mathcal{K}^n$.}

\vspace{0.4cm}

We collect the relations between spherical Crofton measures, generating
functions, and Klain bodies of translation invariant and $\mathrm{SO}(n)$ equivariant
even Minkowski valuations (established in the proofs of
Theorems \ref{main1} and \ref{main2}) \nolinebreak in:

\begin{koro} \label{passing} Let $i \in \{1, \ldots, n - 1\}$ and let $\Phi_i \in \mathbf{MVal}_i^+$ be $\mathrm{SO}(n)$ equivariant and smooth.
If $\mu$ denotes the (smooth) spherical Crofton measure of $\Phi_i$, $\breve{g} \in C^{\infty}(S^{n-1})$ is the generating function of $\Phi_i$,
and $M \in \mathcal{K}^n$ denotes the Klain body of $\Phi_i$, then
\begin{equation} \label{klaincrof17}
h(M,u) = \frac{n\kappa_i\kappa_{n-i}}{2\kappa_{n-1}}{n \choose i}^{-1}\int_{S^{n-1}} |u\cdot v|\,d\mu(v), \qquad u \in S^{n-1},
\end{equation}
and
\[\widehat{\mu}= \frac{2\kappa_{n-1}}{\kappa_i} R_{n-1,i}\,C_{n-1}^{-1}\,\breve{g}.  \]
\end{koro}

\noindent {\bf Examples}:
\begin{enumerate}
\item[(a)] Kiderlen \textbf{\cite{kiderlen05}} proved (in a more general form) that for any $\Phi_1 \in \mathbf{MVal}_1^+$ which is  $\mathrm{SO}(n)$ equivariant and smooth there exists
a unique even zonal function $f \in C^{\infty}(S^{n-1})$ such that for every $K \in \mathcal{K}^n$,
\begin{equation} \label{kidrep}
h(\Phi_1 K,\cdot) = h(K,\cdot) \ast f.
\end{equation}
Using $h(K,u) + h(-K,u) = \mathrm{vol}_1(K|\,\cdot\,)$, we can rewrite (\ref{kidrep}) as
\[h(\Phi_1 K,\cdot) = \mathrm{vol}_1(K|\,\cdot\,) \ast f.  \]
Thus, the function $f$ is (a smooth density of) the spherical Crofton measure of $\Phi_1$.
In order to relate (\ref{kidrep}) with the new representation provided by Theorem \ref{main2},
we recall that (in the sense of distributions)
\begin{equation} \label{s1boxhk}
S_1(K,\cdot) = \Box_n h(K,\cdot),
\end{equation}
where $\Box_n = \frac{1}{n-1}\Delta_S + 1$ and $\Delta_S$ denotes the Laplace--Beltrami operator on $S^{n-1}$ (see, e.g., \textbf{\cite[\textnormal{p.\ 119}]{schneider93}}).
Since $\Box_n$ is a bijection on even functions in $C^{\infty}(S^{n-1})$, there exists an even zonal $\breve{g} \in C^{\infty}(S^{n-1})$ with $\Box_n \breve{g} = f.$
Thus, the $\mathrm{SO}(n)$ equivariance of $\Box_n$ implies (cf.\ \textbf{\cite[\textnormal{p.\ 86}]{grinbergzhang99}}) that
\[h(\Phi_1 K,\cdot) = h(K,\cdot) \ast f = h(K,\cdot) \ast  \Box_n \breve{g} = \Box_n h(K,\cdot) \ast \breve{g} = S_1(K,\cdot) \ast \breve{g}. \]
Finally, the Klain body of $\Phi_1$ can be determined from $f$ by using (\ref{klaincrof17}).

\item[(b)] The first author \textbf{\cite{Schu06a}} proved that $\Phi_{n-1} \in \mathbf{MVal}_{n-1}^+$ is $\mathrm{SO}(n)$ equivariant and smooth
if and only if there exists an $o$-symmetric smooth convex body of revolution $L \in \mathcal{K}^n$ such that for every $K \in \mathcal{K}^n$,
\begin{equation} \label{schurep}
h(\Phi_{n-1} K,\cdot) = S_{n-1}(K,\cdot) \ast h(L,\cdot).
\end{equation}
Consequently, the generating functions of $\mathrm{SO}(n)$ equivariant and smooth Minkowski valuations in $\mathbf{MVal}_{n-1}^+$
are precisely the support functions of $o$-symmetric smooth convex bodies of revolution. It also follows directly from (\ref{schurep}) that
the Klain body of $\Phi_{n-1}$ is given by $2L$.
From Cauchy's projection formula and the fact that the convolution of zonal functions is commutative, it follows that (\ref{schurep}) is equivalent to
\[h(\Phi_{n-1} K,\cdot) = \mathrm{vol}_{n-1}(K|\,\cdot\,) \ast f_L,  \]
where $f_L \in C^{\infty}(S^{n-1})$ is the uniquely determined even zonal function such that
\[h(L,u) = \frac{1}{2}\int_{S^{n-1}}|u\cdot v|f_L(v)\,dv.  \]

\pagebreak

\item[(c)] For $i \in \{1, \ldots, n - 1\}$, let $\Pi_i \in \mathbf{MVal}_i^+$ denote the projection body map of order $i$, defined by
\[h(\Pi_iK,u) = V_i(K|u^{\bot}) = \frac{1}{2} \int_{S^{n-1}} |u\cdot v|\,dS_i(K,v), \qquad u \in S^{n-1}.   \]
Note that each $\Pi_i$ is $\mathrm{SO}(n)$ equivariant but {\it not} smooth. Their (merely) continuous generating function
is given by $g(t) = \frac{1}{2}|t|$. Moreover, it is well known (see, e.g., \textbf{\cite[\textnormal{p.\ 428}]{goodeyweil1}}) that
\[h(\Pi_i K,\cdot) = \frac{\kappa_{n-1}}{\kappa_i} R_{n-i,1}\mathrm{vol}_i^{\bot}(K|\,\cdot\,) =
\frac{\kappa_{n-1}}{\kappa_i}\,\mathrm{vol}_i(K|\,\cdot\,) \ast \widehat{\lambda}_{n-i,1}^{\bot}.  \]
It follows that the Klain body of $\Pi_i$ is a multiple of the Euclidean ball contained in $\bar{e}^{\bot}$, where
$\bar{e} \in \mathrm{Gr}_{i,n}$ is the stabilizer of $\mathrm{S}(\mathrm{O}(i) \times \mathrm{O}(n - i))$.

\item[(d)] For $i \in \{2, \ldots, n\}$, let $\mathrm{M}_i \in \mathbf{MVal}_{n+1-i}$ denote the {\it normalized} mean section operator of order $i$,
introduced by Goodey and Weil \textbf{\cite{goodeyweil1,
goodeyweil2}} and given by
\begin{equation} \label{MiK}
h(\mathrm{M}_iK,\cdot) = \int_{\mathrm{AGr}_{i,n}} h(\mathrm{J}(K \cap E),\cdot)\,d\mu_i(E).
\end{equation}
Here $\mathrm{J} \in \mathbf{MVal}_1$ is defined by $\mathrm{J}K =
K - s(K)$, where $s: \mathcal{K}^n \rightarrow \mathbb{R}^n$ is
the Steiner point map (see, e.g., \textbf{\cite[\textnormal{p.\
50}]{schneider93}}) and $\mu_i$ is the invariant measure on $\mathrm{AGr}_{i,n}$ normalized such that the set of
planes having \linebreak non-empty intersection with the Euclidean unit ball has measure $\kappa_{n-i}$.
Note that $\mathrm{M}_i$ is \emph{not} even. However, it was
proved in \textbf{\cite{goodeyweil1}} that
\[h(\mathrm{M}_iK,\cdot) + h(\mathrm{M}_i(-K),\cdot)=\frac{i\kappa_i\kappa_{n-1}}{n\kappa_{i-1}\kappa_n}\,R_{n+1-i,1}\mathrm{vol}_{n+1-i}(K|\,\cdot\,).  \]
Thus, a multiple of $\widehat{\lambda}_{n+1-i,1}$ is the spherical Crofton measure for the even part of $\mathrm{M}_i$ and its Klain body is
a multiple of the Euclidean ball contained in the subspace $\bar{e} \in \mathrm{Gr}_{n+1-i,n}$.

Goodey and Weil \textbf{\cite{goodeyweil2}} also determined the family of
generating functions for the mean section operators.
In order to explain their result, recall that in Berg's solution of the Christoffel--Minkowski problem (see, e.g.,
\textbf{\cite{goodeyyaskin, schneider93}}) he proved the following: For every $n \geq 2$ there exists a uniquely
determined $C^{\infty}$ function $\zeta_n$ on $(-1,1)$ such that the associated zonal function $\breve{\zeta}_n \in L^1(S^{n-1})$ is orthogonal
to the restriction of all linear functions to $S^{n-1}$ and satisfies, for every $K \in \mathcal{K}^n$,
\begin{equation} \label{bergfct}
h(\mathrm{J}K,\cdot) =  S_1(K,\cdot) \ast \breve{\zeta}_n.
\end{equation}

\pagebreak

Goodey and Weil \textbf{\cite[\textnormal{Theorem 4.4}]{goodeyweil2}} proved that
\begin{equation} \label{genfctmi}
h(\mathrm{M}_iK,\cdot) = q_{n,i}\,S_{n+1-i}(K,\cdot) \ast
\breve{\zeta}_i,
\end{equation}
where
\begin{equation} \label{qni}
q_{n,i} = \frac{i-1}{2\pi(n+1-i)}\,\frac{\kappa_{i-1}\kappa_{i-2}\kappa_{n-i}}{\kappa_{i-3}\kappa_{n-2}}.
\end{equation}
\end{enumerate}

Let $C_\mathrm{o}^{\infty}(S^{n-1})$ denote the subspace of smooth functions on $S^{n-1}$ which are orthogonal to the restriction of all linear functions.
It is well known that the linear differential operator $\Box_n: C_\mathrm{o}^{\infty}(S^{n-1}) \rightarrow C_\mathrm{o}^{\infty}(S^{n-1})$ is an isomorphism (see, e.g., \textbf{\cite{groemer96}}).
Moreover, since every twice continuously differentiable function on $S^{n-1}$ is a difference of support functions (see, e.g.,
\textbf{\cite[\textnormal{p.\ 49}]{schneider93}}), it follows from (\ref{s1boxhk}) and (\ref{bergfct}) that
\begin{equation} \label{boxninverse}
f = (\Box_n f) \ast \breve{\zeta}_n
\end{equation}
for every $f \in C_\mathrm{o}^{\infty}(S^{n-1})$. In the next section we need the following more general fact.

\begin{prop} \label{propboxj} For every $n \geq 2$ and $2 \leq j \leq n$, the integral transform
\[\mathrm{F}_{\zeta_j}: C_{\mathrm{o}}^{\infty}(S^{n-1}) \rightarrow C_{\mathrm{o}}^{\infty}(S^{n-1}), \quad f \mapsto f \ast \breve{\zeta}_j,  \]
is an isomorphism.

\end{prop}

Proposition \ref{propboxj} can be deduced from a result of Goodey and Weil \linebreak \textbf{\cite[\textnormal{Theorem 4.3}]{goodeyweil2}}.
A different and more elementary proof was given very recently in \textbf{\cite{BPSW2014}}. Proposition \ref{propboxj} and (\ref{boxninverse}) give rise to the following:

\vspace{0.3cm}

\noindent {\bf Definition} \emph{For $2 \leq j \leq n$, let $\Box_j: C_{\mathrm{o}}^{\infty}(S^{n-1}) \rightarrow
C_{\mathrm{o}}^{\infty}(S^{n-1}) $ denote the linear operator which is
inverse to the integral transform $\mathrm{F}_{\zeta_j}$.}

\vspace{1cm}

\centerline{\large{\bf{ \setcounter{abschnitt}{6}
\arabic{abschnitt}. The Hard Lefschetz operators}}}

\reseteqn \alpheqn \setcounter{theorem}{0}

\vspace{0.6cm}

In this final section we determine the action of Alesker's Hard Lefschetz \linebreak operators on translation invariant and
$\mathrm{SO}(n)$ equivariant even Minkowski valuations in terms
of corresponding integral transforms of the generating function, the spherical Crofton measure, and the support function of the Klain body of the respective Minkowski valuation.
Our investigations are motivated by recent applications of the derivation operator on Minkowski valuations in the theory of geometric inequalities (see \textbf{\cite{parapschu}}).

\pagebreak

In the recent article \textbf{\cite{parapschu}}, Parapatits and the first author  showed that for any $\Phi \in \mathbf{MVal}^{(+)}$ there exist
$\Phi^{(j)} \in \mathbf{MVal}^{(+)}$, where $0 \leq j \leq n$, such that
\[\Phi(K + rB) = \sum_{j=0}^n r^{n-j}\Phi^{(j)}(K)  \]
for every $K \in \mathcal{K}^n$ and $r \geq 0$.
This Steiner-type formula, in turn, gives rise to the definition of a derivation operator $\Lambda: \mathbf{MVal} \rightarrow \mathbf{MVal}$.

\vspace{0.3cm}

\noindent {\bf Definition} \emph{For $\Phi \in \mathbf{MVal}$, define $\Lambda \Phi \in \mathbf{MVal}$ by}
\[ h((\Lambda \Phi)(K),u) = \left . \frac{d}{dt} \right |_{t=0} h(\Phi(K + tB),u), \qquad u \in S^{n-1}.  \]

\vspace{0.2cm}

Note that $\Lambda$ commutes with the action of $\mathrm{SO}(n)$ and that $\Lambda$ preserves parity. Moreover, if $\Phi_i \in \mathbf{MVal}_i^{(+)}$ is $\mathrm{SO}(n)$ equivariant and smooth, then
so is $\Lambda\Phi_i \in \mathbf{MVal}_{i-1}^{(+)}$. We also emphasize that if $\varphi_i \in \mathbf{Val}_i^{(+)}$ is the real valued valuation associated with $\Phi_i$,
then $\Lambda \varphi_i \in \mathbf{Val}_{i-1}^{(+)}$ is associated with $\Lambda\Phi_i$.

\vspace{0.3cm}

\noindent {\bf Example}:

\vspace{0.1cm}

\noindent By definition of the projection body maps of order $i$ and (\ref{steiner}), we have
\[\Lambda^{n-1-i}\Pi_{n-1}=\frac{(n-1)!}{i!}\,\Pi_i.\]

\vspace{0.1cm}

\begin{theorem} \label{der17} Suppose that $i \in \{2, \ldots, n - 1\}$ and let $\Phi_i \in \mathbf{MVal}_i^+$ be \linebreak $\mathrm{SO}(n)$ equivariant and smooth.
\begin{enumerate}
\item[(i)] If $\breve{g} \in C^{\infty}(S^{n-1})$ is the generating function of $\Phi_i$, then the generating function of $\Lambda \Phi_i$ is given by $i\,\breve{g}$.
\item[(ii)] If $\mu$ is the (smooth) spherical Crofton measure of $\Phi_i$, then the spherical Crofton measure $\nu$ of $\Lambda \Phi_i$ is determined by
\[\widehat{\nu} = \frac{i\kappa_i}{\kappa_{i-1}}\, R_{i,i-1}\,\widehat{\mu}.  \]
\item[(iii)] If $M \in \mathcal{K}^n$ is the Klain body of $\Phi_i$, then the Klain body $N \in \mathcal{K}^n$ of $\Lambda \Phi_i$ is determined by
\[\widehat{h(N,\cdot)} = \frac{(n-i+1)\kappa_{n-i+1}}{\kappa_{n-i}}\, R_{i,i-1}\,\widehat{h(M,\cdot)}.   \]
\end{enumerate}
\end{theorem}

\noindent {\it Proof.} Since (\ref{steiner}) can be generalized to arbitrary area measures of order $i$, more precisely,
\[S_i(K+tB,\cdot) = \sum_{j=0}^i t^{i-j} {i \choose j} S_j(K,\cdot),   \]
we have
\[ h(\Phi_i(K + tB,\cdot))=S_i(K+tB,\cdot) \ast \breve{g} = \sum_{j=0}^i t^{i-j} {i \choose j} S_j(K,\cdot)\ast \breve{g}  \]
and, thus, $h((\Lambda \Phi_i)K,\cdot) = S_{i-1}(K,\cdot) \ast i\, \breve{g}$, which proves (i).

\vspace{0.1cm}

By Corollary \ref{passing}, on one hand
\[\breve{g} = \frac{\kappa_i}{2\kappa_{n-1}}\, C_{n-1\,}R_{n-1,i}^{-1}\, \widehat{\mu}   \]
and, by (i) and Corollary \ref{passing} again, on the other hand
\[\widehat{\nu} = \frac{2i\kappa_{n-1}}{\kappa_{i-1}}\, R_{n-1,i-1\,}C_{n-1}^{-1}\,\breve{g}  \]
which together with $R_{n-1,i-1} = R_{i,i-1} \circ R_{n-1,i}$ prove claim (ii).

\vspace{0.1cm}

In order to prove (iii), we use (\ref{proof2}) and (ii) to arrive at
\[\widehat{h(N,\cdot)} = C_{i-1}\widehat{\nu} = \frac{i\kappa_i}{\kappa_{i-1}}\, C_{i-1\,}R_{i,i-1}\,\widehat{\mu}.   \]
Thus, from another application of (\ref{proof2}), we obtain
\[ \widehat{h(N,\cdot)} = \frac{i\kappa_i}{\kappa_{i-1}}\, C_{i-1\,}R_{i,i-1\,}C_{i}^{-1} \widehat{h(M,\cdot)}.  \]
Using now (\ref{rijci}) gives the desired result. \hfill $\blacksquare$

\vspace{0.4cm}

Theorem \ref{der17} yields a set of necessary conditions for a function or a measure to be the generating function or the spherical Crofton measure, respectively, of a Minkowski valuation.
For example, for the generating function $\breve{g}$ of a Minkowski valuation $\Phi_i \in \mathbf{MVal}_i$, all the functions $S_j(K,\cdot) \ast \breve{g}$, $1 \leq j \leq i$,
must be support functions for every $K \in \mathcal{K}^n$.

\vspace{0.2cm}

Next, we turn to the integration operator. Formula (\ref{hardlefint}) of Bernig \textbf{\cite{Bernig07}} (which we used as a definition of $\mathfrak{L}$ on $\mathbf{Val}$) motivates the following:

\noindent {\bf Definition} \emph{For $\Phi \in \mathbf{MVal}$, define $\mathfrak{L} \Phi \in \mathbf{MVal}$ by}
\[ h((\mathfrak{L} \Phi)(K),u) = \int_{\mathrm{AGr}_{n-1,n}} h(\Phi(K \cap E),u)\,d\sigma_{n-1}(E), \qquad u \in S^{n-1}.  \]

\vspace{0.3cm}

Note that $\mathfrak{L}$ commutes with the action of $\mathrm{SO}(n)$ and that $\mathfrak{L}$ preserves parity.
Moreover, if $\Phi_i \in \mathbf{MVal}_i^{(+)}$ is $\mathrm{SO}(n)$ equivariant and smooth, then so is
$\mathfrak{L}\Phi_i \in \mathbf{MVal}_{i+1}^{(+)}$. We also emphasize that if $\varphi_i \in \mathbf{Val}_i^{(+)}$ is the
real valued valuation associated with $\Phi_i$, then $\mathfrak{L}\varphi_i \in \mathbf{Val}_{i+1}^{(+)}$ is associated with $\mathfrak{L}\Phi_i$.

\vspace{0.4cm}

\noindent {\bf Example}:

\vspace{0.2cm}

\noindent For $k \in \{1, \ldots, n\}$, it follows by induction on $k$ from a well known formula of Crofton (see, e.g., \textbf{\cite[\textnormal{p.\ 124}]{Klain:Rota}}) that
\[\int_{\mathrm{AGr}_{n-1,n}} \!\!\!\!\!\!\!\!\! \cdots \int_{\mathrm{AGr}_{n-1,n}}\!\!\!\!\!\!\!\!\!\!\!\!\!\!\!\! f(E_1 \cap \cdots \cap E_k)\,d\sigma_{n-1}(E_1) \cdots d\sigma_{n-1}(E_k)
= \frac{k!\kappa_k}{2^k}\!\! \int_{\mathrm{AGr}_{n-k,n}} \!\!\!\!\!\!\!\!\!\!\!\!\!\!\!\! f(F)\,d\sigma_{n-k}(F)   \]
for every $f \in L^1(\mathrm{AGr}_{n-k,n})$. Consequently, for $\Phi \in \mathbf{MVal}$,
\[h((\mathfrak{L}^k\Phi)(K),u) = \frac{k!\kappa_k}{2^k} \int_{\mathrm{AGr}_{n-k,n}} h(\Phi(K \cap F),u)\,d\sigma_{n-k}(F), \qquad u \in S^{n-1}.   \]
By definition (\ref{MiK}) of the mean section operator of order $i$ and the fact that
\[\sigma_{n-i} = \left[\begin{array}{c} n\\i \end{array}\right] \mu_{n-i},  \]
we conclude that for $i \in \{0, \ldots, n - 2\}$,
\begin{equation} \label{LiJMn-i}
\mathfrak{L}^i\mathrm{J} = \frac{i!\kappa_i}{2^i} \left[\begin{array}{c} n\\i \end{array}\right] \mathrm{M}_{n-i} = \frac{n!\kappa_n}{2^i(n-i)!\kappa_{n-i}} \mathrm{M}_{n-i}.
\end{equation}

\vspace{0.1cm}

\begin{theorem} \label{int17} Suppose that $i \in \{1, \ldots, n - 2\}$ and let $\Phi_i \in \mathbf{MVal}_i^+$ be \linebreak $\mathrm{SO}(n)$ equivariant and smooth.
\begin{enumerate}
\item[(i)] If $\breve{g} \in C^{\infty}(S^{n-1})$ is the generating function of $\Phi_i$, then the generating function $\breve{f} \in C^{\infty}(S^{n-1})$ of $\mathfrak{L} \Phi_i$ is given by
\[\breve{f} = \frac{(n-i)\kappa_{i+1}\kappa_{n-i}}{2\kappa_i\kappa_{n-i-1}}\,R_{n-1,i+1\,}^{-1}R_{i,i+1\,}R_{n-1,i}\,\breve{g}.   \]
\item[(ii)] If $\mu$ is the (smooth) spherical Crofton measure of $\Phi_i$, then the spherical Crofton measure $\nu$ of $\mathfrak{L} \Phi_i$ is determined by
\[\widehat{\nu} = \frac{(n-i)\kappa_{n-i}}{2\kappa_{n-i-1}}\, R_{i,i+1}\,\widehat{\mu}.  \]
\item[(iii)] If $M \in \mathcal{K}^n$ is the Klain body of $\Phi_i$, then the Klain body $N \in \mathcal{K}^n$ of $\mathfrak{L} \Phi_i$ is determined by
\[\widehat{h(N,\cdot)} = \frac{(i+1)\kappa_{i+1}}{2\kappa_i}\, R_{i,i+1}\,\widehat{h(M,\cdot)}.   \]
\end{enumerate}
\end{theorem}
{\it Proof.} We begin with the proof of (iii). To this end let $\varphi_i \in \mathbf{Val}_i^+$ denote the associated real valued valuation of $\Phi_i$.
First note that, by (\ref{hardfourier}),
\begin{equation} \label{lfourlam}
\mathfrak{L}\varphi_i=\frac{1}{2}\mathbb{F}\Lambda\mathbb{F}\varphi_i
\end{equation}
is the associated real valued valuation of $\mathfrak{L}\Phi_i$. Moreover, by (\ref{proof2}),
\[\mathrm{Kl}_i\varphi_i=\widehat{h(M,\cdot)} \qquad \mbox{and} \qquad \mathrm{Kl}_{i+1}(\mathfrak{L}\varphi_i)=\widehat{h(N,\cdot)}.\]

Thus, (\ref{lfourlam}), (\ref{deffourier}), Theorem \ref{der17}
(iii), and (\ref{rijbot}) yield
\[\widehat{h(N,\cdot)}=\frac{1}{2} \left (\mathrm{Kl}_{n-i-1}(\Lambda \mathbb{F}\varphi_i) \right )^{\bot} =\frac{(i+1)\kappa_{i+1}}{2\kappa_i} R_{i,i+1}\,\widehat{h(M,\cdot)}. \]

\vspace{0.1cm}

In order to prove (ii), we use Corollary \ref{passing}, Lemma \ref{keylem} (i), and part (iii) which we just proved, to obtain
\[\widehat{\nu} =  \frac{(i+1)\kappa_{i+1}}{2\kappa_i}\,C_{i+1\,}^{-1}R_{i,i+1\,}C_i\, \widehat{\mu}. \]
An application of (\ref{rijci}) gives the desired result.

\vspace{0.1cm}

Finally, by Corollary \ref{passing}, we have
\[\widehat{\mu} = \frac{\kappa_{n-1}}{\kappa_i}\, R_{n-1,i\,}C_{n-1}^{-1}\,\breve{g} \qquad \mbox{and} \qquad  \breve{f}=\frac{\kappa_{i+1}}{\kappa_{n-1}}\,C_{n-1\,}R_{n-1,i+1}^{-1} \widehat{\nu}.  \]
Thus, using (ii), we arrive at
\[\breve{f} = \frac{(n-i)\kappa_{i+1}\kappa_{n-i}}{2\kappa_i\kappa_{n-i-1}}\,C_{n-1\,}R_{n-1,i+1\,}^{-1}R_{i,i+1\,} R_{n-1,i\,}C_{n-1}^{-1}\,\breve{g}.   \]
Using now (\ref{rijci}) three times gives (i).
\hfill $\blacksquare$

\pagebreak

We note that it is an open problem whether the Alesker--Fourier
transform $\mathbb{F}$ is well defined for even Minkowski
valuations in $\mathbf{MVal}$ that are $\mathrm{SO}(n)$
equivariant and smooth. More precisely, if $\mu$ is the spherical Crofton measure of a smooth and $\mathrm{SO}(n)$
equivariant Minkowski valuation $\Phi_i \in \mathbf{MVal}_{i}^{+}$, then it is not known in general whether $\mu^{\bot}$
is the spherical Crofton measure of a smooth and $\mathrm{SO}(n)$
equivariant Minkowski valuation $\mathbb{F}\Phi_i \in \mathbf{MVal}_{n-i}^{+}$.

\vspace{0.2cm}

In the last part of this section, we use  Proposition \ref{propboxj} and relation (\ref{LiJMn-i}) to deduce a
more explicit expression for the generating function of
$\mathfrak{L}\Phi_i$ in terms of the generating function of an
$\mathrm{SO}(n)$ equivariant smooth $\Phi_i \in \mathbf{MVal}_i$. Since we
want to prove this result for Minkowski valuations which are {\it not necessarily even},
recall that a zonal function $\breve{g} \in C^{\infty}_{\mathrm{o}}(S^{n-1})$ is called a \emph{generating function} for $\Phi_i$ if
\[h(\Phi_i K,\cdot) = S_i(K,\cdot) \ast \breve{g}  \]
holds for every $K \in \mathcal{K}^n$. Note that if $\Phi_i \in \mathbf{MVal}_{i}$ admits a generating function $\breve{g} \in C^{\infty}_{\mathrm{o}}(S^{n-1})$, then it follows from well known density properties of area measures of convex bodies that $\breve{g}$ is uniquely determined by $\Phi_i$. Moreover, by Theorem \ref{main2}, every smooth $\mathrm{SO}(n)$ equivariant \emph{even} $\Phi_i \in \mathbf{MVal}_{i}$ admits a generating function.

\begin{theorem} Suppose that $i \in \{1, \ldots, n - 2\}$ and let $\Phi_i \in \mathbf{MVal}_i$ be $\mathrm{SO}(n)$ equivariant and smooth.
If $\Phi_i$ admits a generating function $\breve{g} \in C^{\infty}(S^{n-1})$, then $\mathfrak{L} \Phi_i$ also admits a generating function $\breve{f} \in
C^{\infty}(S^{n-1})$ which is given by
\[\breve{f} = c_{n,i\,} \Box_{n-i+1} \breve{g} \ast \breve{\zeta}_{n-i},   \]
where
\[c_{n,i} = \frac{i(n-i-1)(n-i+1)\kappa_{n-i-2}^2 \kappa_{n-i+1}\kappa_i}{2(n-i)(i+1)\kappa_{n-i-3}\kappa_{n-i}^2 \kappa_{i-1}}.  \]
\end{theorem}
{\it Proof.} Consider the subspace of $\mathbf{Val}_i^{\infty}$
spanned by valuations of the form
\[\psi_i(K)= \int_{S^{n-1}} h(u)\,dS_i(K,u), \qquad K \in \mathcal{K}^n,   \]
where $h \in C_{\mathrm{o}}^{\infty}(S^{n-1})$. Note that $h$ is uniquely determined by $\psi_i$.
By the $\mathrm{SO}(n)$ equivariance and linearity of $\mathfrak{L}: \mathbf{Val}_i^{\infty} \rightarrow \mathbf{Val}_{i+1}^{\infty}$, it follows that
there exists a linear operator $\mathrm{T}_i: C^{\infty}_{\mathrm{o}}(S^{n-1})
\rightarrow C^{\infty}_{\mathrm{o}}(S^{n-1})$ which is $\mathrm{SO}(n)$
equivariant, such that
\[(\mathfrak{L}\psi_i)(K) = \int_{S^{n-1}}(\mathrm{T}_ih)(u)\, \,dS_{i+1}(K,u), \qquad K \in \mathcal{K}^n. \]

\pagebreak

If $\varphi_i \in \mathbf{Val}_i^{\infty}$ denotes the
real valued associated valuation of $\Phi_i$, then
\[\varphi_i(K) =  \int_{S^{n-1}} \breve{g}(u)\,dS_i(K,u), \qquad K \in \mathcal{K}^n. \]
Therefore, we have to show that for every $h \in C_{\mathrm{o}}^{\infty}(S^{n-1})$,
\begin{equation} \label{toshow}
\mathrm{T}_ih = c_{n,i\,} \Box_{n-i+1} h \ast \breve{\zeta}_{n-i} = c_{n,i\,} \mathrm{F}_{\zeta_{n-i}}(\Box_{n-i+1} h).
\end{equation}

To this end, let $h \sim \sum_{k\geq 0}H_k$ be the series expansion of $h$ into spherical harmonics.
By the linearity and $\mathrm{SO}(n)$ equivariance of $\mathrm{T}_i$,
there exists a sequence $a_k^n[T_i] \in \mathbb{R}$ of multipliers such
that the spherical harmonics expansion of $\mathrm{T}_ih$ is given
by $\mathrm{T}_ih \sim \sum_{k\geq 0}a_k^n[T_i]\,H_k$. In particular, the operator $\mathrm{T}_i$ is
uniquely determined by the sequence $a_k^n[T_i]$. Moreover, the sequence $a_k^n[T_i]$ is determined by the $\mathrm{T}_i$-image of any function being the sum of nonzero
harmonics of all orders different from 1. Hence, since, by the Funk--Hecke Theorem (see, e.g., \textbf{\cite[\textnormal{p.\ 98}]{groemer96}}), $\mathrm{F}_{\zeta_{n-i}}$ and $\Box_{n-i+1}$ are also multiplier transformations, it suffices to prove (\ref{toshow}) for one such function.

By Proposition \ref{propboxj}, each of Berg's functions $\breve{\zeta}_j$, $2 \leq j \leq n$,
is such a sum of nonzero harmonics of all orders different from 1. However, $\breve{\zeta}_j$ is not smooth but merely in $L^1_{\mathrm{o}}(S^{n-1})$, the subspace
of $L^1(S^{n-1})$ consisting of functions which are orthogonal to all spherical harmonics of degree 1. This is not a problem because,
as multiplier transformations, $\mathrm{T}_i$ as well as $\Box_j$ are selfadjoint and thus can be extended in the sense of distributions to $L^1_{\mathrm{o}}(S^{n-1})$.
Therefore, on one hand
\begin{equation} \label{boxjid}
f= (\Box_j f) \ast \breve{\zeta}_j
\end{equation}
holds in the sense of distributions for every $f \in L^1_{\mathrm{o}}(S^{n-1})$ and $j \in \{2, \ldots, n\}$. On the other hand, since, by (\ref{LiJMn-i}),
\[\mathfrak{L}\mathrm{M}_{n+1-i} = \frac{(n-i+1)\kappa_{n-i+1}}{2\kappa_{n-i}}\,\mathrm{M}_{n-i},  \]
we obtain from (\ref{genfctmi}) that in the sense of distributions
\[\mathrm{T}_i(q_{n,n-i+1}\breve{\zeta}_{n-i+1}) = \frac{(n-i+1)\kappa_{n-i+1}}{2\kappa_{n-i}}\,q_{n,n-i\,} \breve{\zeta}_{n-i}.  \]
Using (\ref{boxjid}) with $f = \breve{\zeta}_{n-i}$ and $j = n - i+1$ and that multiplier transformations and convolutions of zonal functions are commutative, this is equivalent to
\[\mathrm{T}_i\breve{\zeta}_{n-i+1} = \frac{(n-i+1)\kappa_{n-i+1}q_{n,n-i}}{2\kappa_{n-i}q_{n,n-i+1}} (\Box_{n-i+1} \breve{\zeta}_{n-i+1}) \ast \breve{\zeta}_{n-i}.    \]
Plugging in the explicit values of $q_{n,m}$ given in (\ref{qni}) completes the proof. \hfill
$\blacksquare$

\pagebreak

\centerline{\large{\bf{\setcounter{abschnitt}{7} Appendix}}}

\renewcommand{\theequation}{A.\arabic{equation}}{\setcounter{equation}{0}
\setcounter{theorem}{0}
\renewcommand{\thesection}{\Alph{section}}
\setcounter{section}{1}

\vspace{0.6cm}

The purpose of this appendix is the proof of the proposition below, which (together with Theorem 2) shows that the
notion of spherical Crofton measure of a Minkowski valuation is independent from the identification of $S^{n-1}$ with $\mathrm{O}(n)/\mathrm{O}(n - 1)$
or $\mathrm{SO}(n)/\mathrm{SO}(n - 1)$, respectively.

\vspace{0.3cm}

\noindent {\bf Proposition A} \emph{If $f \in C(S^{n-1})$ is $\mathrm{S}(\mathrm{O}(i) \times \mathrm{O}(n - i))$ invariant, then $f$ is also
$\mathrm{O}(i) \times \mathrm{O}(n - i)$ invariant.   }

\vspace{0.2cm}

\noindent {\it Proof}. Let $H_1=\mathrm{O}(i) \times \mathrm{O}(n - i)$, $H_2=\mathrm{S}(\mathrm{O}(i) \times \mathrm{O}(n - i))$ and let
$\Gamma$ be an arbitrary $\mathrm{O}(n)$ irreducible subspace of $L^2(S^{n-1})$.
Recall that we denote by $\Gamma^{H_i}$ the subspace of $H_i$ invariant
elements of $\Gamma$. Since $L^2(S^{n-1})$ is an orthogonal direct sum of $\mathrm{O}(n)$ irreducible subspaces,
it will be sufficient to prove that $\Gamma^{H_1}=\Gamma^{H_2}$.

Since $H_2 \subset H_1$, we obviously have $ \Gamma^{H_1} \subset \Gamma^{H_2}$.
Moreover, it follows from the Frobenius reciprocity theorem (see, e.g., \textbf{\cite[\textnormal{Theorem 9.9}]{knapp}}) that
\[\dim \Gamma^{H_1}= \dim \mathrm{Hom}_{\mathrm{O}(n)}(\Gamma,C(\mathrm{Gr}_{i,n}))\]
and
\begin{equation}\label{eq:V_K}
\dim \Gamma^{H_2} =  \mathrm{Hom}_{\mathrm{SO}(n)}(\mathrm{Res}^{\mathrm{O}(n)}_{\mathrm{SO}(n)} \Gamma, C(\mathrm{Gr}_{i,n})).
\end{equation}
Here, $\mathrm{Hom}_G$ denotes the space of linear $G$-equivariant maps and $\mathrm{Res}^{\mathrm{O}(n)}_{\mathrm{SO}(n)}$ denotes the restriction
of an $\mathrm{O}(n)$ representation to $\mathrm{SO}(n)$.

From the description of the irreducible representations of $\mathrm{O}(n)$
in terms of the irreducible representations of $\mathrm{SO}(n)$ (see, e.g., \textbf{\cite[\textnormal{Lemma 3.1}]{ABS2011}})
and \eqref{decompsn1}, we obtain
\[\mathrm{Res}^{\mathrm{O}(n)}_{\mathrm{SO}(n)} \Gamma = \Gamma_{(k,0,\ldots,0)}\]
for some $k \in \mathbb{N}$. Thus, by \eqref{eq:V_K} and (\ref{decompgri}), we have $\dim \Gamma^{H_2}=0$ if $k$ is odd and $\dim \Gamma^{H_2}=1$ if $k$ is even.
Hence, the proposition will be proved if we can show that
\[\dim \mathrm{Hom}_{\mathrm{O}(n)}(\Gamma, C(\mathrm{Gr}_{i,n})) \geq 1\]
whenever $\mathrm{Res}^{\mathrm{O}(n)}_{\mathrm{SO}(n)} \Gamma = \Gamma_{(2k,0,\ldots,0)}$.
But if we identify even functions on $S^{n-1}$ with functions on $\mathrm{Gr}_{1,n}$, then the restriction of the
Radon transform $R_{1,i}: C(\mathrm{Gr}_{1,n}) \to C(\mathrm{Gr}_{i,n})$ is a non-trivial $\mathrm{O}(n)$ equivariant map from
$\Gamma$ to $C(\mathrm{Gr}_{i,n})$. \hfill $\blacksquare$

\vspace{1cm}

\noindent {{\bf Acknowledgments} The work of the authors was
supported by the Austrian \linebreak Science Fund (FWF), Project number:
P\,22388-N13. The first author was also supported by the European Research Council (ERC), Project number: 306445.
The second author was also supported by the German Research Foundation (DFG), Project number: BE 2484/5-1.

\begin{small}

\[ \begin{array}{ll} \mbox{Franz Schuster} & \mbox{Thomas Wannerer} \\
\mbox{Vienna University of Technology \phantom{www}} & \mbox{Goethe-University Frankfurt} \\ \mbox{franz.schuster@tuwien.ac.at} & \mbox{wannerer@mathematik.uni-frankfurt.de}
\end{array}\]

\end{small}

\end{document}